\documentclass[11pt]{amsart}

\usepackage{amssymb}

\usepackage{upref}

\newcommand{\conservepaper}{
 \hoffset=-0.75in
 \setlength{\textwidth}{6.5in}
 \voffset=-0.5in
 \setlength{\textheight}{10.0in}
 \setlength{\textheight}{9.0in} 
 }
 \conservepaper

\swapnumbers 

\numberwithin{equation}{section}
\newtheorem{thm}[equation]{Theorem}

\newtheorem{lem}[equation]{Lemma}
\newtheorem{cor}[equation]{Corollary}
\newtheorem{prop}[equation]{Proposition}

\theoremstyle{definition}

\newtheorem{defn}[equation]{Definition}
\newtheorem{conj}[equation]{Conjecture}

\newtheorem{assump}[equation]{Assumption}
\newtheorem{notation}[equation]{Notation}

\theoremstyle{remark}

\newtheorem{rem}[equation]{Remark}        
         
\newtheorem{claim}[equation]{Claim}
\newtheorem{ack}[equation]{Acknowledgment}  

\newcommand{\pref}[1]{{\upshape(}\ref{#1}{\upshape)}}
\renewcommand{\see}[1]{{\upshape(}see~\ref{#1}{\upshape)}}
\newcommand{\cf}[1]{(cf.~\ref{#1})}
\newcommand{\fullref}[2]{\ref{#1}\pref{#1#2}}

\newcommand{\Lie}[1]{\mathfrak{#1}}
\newcommand{\diag}{\operatorname{diag}}

\newcommand{\SL}{\operatorname{SL}}
\newcommand{\GL}{\operatorname{GL}}
\newcommand{\PGL}{\operatorname{PGL}}
\newcommand{\SO}{\operatorname{SO}}
\newcommand{\SU}{\operatorname{SU}}
\newcommand{\On}{\operatorname{O}}
\renewcommand{\Sp}{\operatorname{Sp}}
\newcommand{\so}{\operatorname{\Lie{so}}}
\newcommand{\sltwoC}{\operatorname{\Lie{sl}}(2,\complex)}
\newcommand{\sltwoR}{\operatorname{\Lie{sl}}(2,\real)}
\newcommand{\Id}{\operatorname{Id}}

\newcommand{\Rad}{\operatorname{Rad}}
\newcommand{\real}{\mathord{\mathbb{R}}}
\newcommand{\complex}{\mathord{\mathbb{C}}}

\newcommand{\semiprod}{\ltimes}
\newcommand{\closure}[1]{\overline{#1}}
\newcommand{\iso}{\cong}

\newcommand{\Rrank}{\mathop{\real\text{\upshape-rank}}}

\newcommand{\GamPrime}{\mathord{\Gamma'}}
\newcommand{\gamprime}{\mathord{\gamma'}}
\newcommand{\muLZ}{\mu_{LZ}}
\newcommand{\mul}{\mu_L}
\newcommand{\G}{\mathord{\mathbb G}}

\newcommand{\muH}[2]{\bigl[ #1, #2 \bigr]}

\newcommand{\bigset}[2]{\left\{\, #1 
 \mathrel{\left| \vphantom {\left\{ #1 \mid #2 \right\} } \right.}
 #2 \,\right\} }

 \newcounter{case}
 \newenvironment{case}[1][\unskip]{\refstepcounter{case}\em
 \medskip \noindent Case \thecase\ #1.\ }{\unskip\upshape}
 \renewcommand{\thecase}{\arabic{case}}

 \newcounter{subcase}
 \newenvironment{subcase}[1][\unskip]{\refstepcounter{subcase}\em
 \medskip \noindent Subcase \thesubcase\ #1.\ }{\unskip\upshape}
\numberwithin{subcase}{case}

\begin{document}

\title[Clifford-Klein forms of $\SO(2,n)/H$]
 {Compact Clifford-Klein forms 
 \\ of homogeneous spaces of $\SO(2,n)$
 }

\author{Hee Oh}
 \address{Department of Mathematics, Oklahoma State University,
Stillwater, OK 74078}
 \curraddr{Institute of Mathematics, The Hebrew University,
 Jerusalem 91904, Israel}
 \email{heeoh@math.huji.ac.il}

\author{Dave Witte}
 \address{Department of Mathematics, Oklahoma State University,
Stillwater, OK 74078}
 \email{dwitte@math.okstate.edu,
  http://www.math.okstate.edu/\char'176dwitte} 


\date{March 3, 1999 \bf (Corrected version)} 

\begin{abstract}
 A homogeneous space $G/H$ is said to have a compact Clifford-Klein form if
there exists a discrete subgroup~$\Gamma$ of~$G$ that acts properly
discontinuously on $G/H$, such that the quotient space $\Gamma \backslash
G/H$ is compact. When $n$ is even, we find every closed, connected
subgroup~$H$ of $G = \SO(2,n)$, such that $G/H$ has a compact Clifford-Klein
form, but our classification is not quite complete when $n$~is odd.  The
work reveals new examples of homogeneous spaces of $\SO(2,n)$ that have
compact Clifford-Klein forms, if $n$ is even. Furthermore, we show
that if $H$ is a closed, connected subgroup of $G = \SL(3,\real)$, and
neither  $H$ nor $G/H$ is compact, then $G/H$ does not have a compact
Clifford-Klein form, and we also study noncompact Clifford-Klein forms of
finite volume.
 \end{abstract}

\maketitle

\section{Introduction}

\begin{assump}
 Throughout this paper, $G$ is a Zariski-connected, almost simple, linear Lie
group (``Almost simple" means that every proper normal subgroup of~$G$
either is finite or has finite index.) In the main results, $G$ is assumed
to be $\SO(2,n)$ (with $n \ge 3$).  

There would be no essential loss of generality if one were to require $G$ to
be connected, instead of only Zariski connected \see{Gconn}. However,
$\SO(2,n)$ is not connected (it has two components) and the authors prefer
to state results for $\SO(2,n)$, instead of for the identity component of
$\SO(2,n)$.
 \end{assump}

\begin{defn}
 Let $H$ be a closed, connected subgroup of~$G$. We say that the
homogeneous space $G/H$ has a \emph{compact Clifford-Klein form} if there
is a discrete subgroup~$\Gamma$ of~$G$, such that 
 \begin{itemize}
 \item $\Gamma$ acts properly on~$G/H$; and
 \item $\Gamma \backslash G/H$ is compact.
 \end{itemize}
 (Alternative terminology would be to say that $G/H$ has a
\emph{tessellation}, because the $\Gamma$-translates of a fundamental
domain for $\Gamma \backslash G/H$ tessellate $G/H$, or one could simply
say that $G/H$ has a \emph{compact quotient}.)
 See the surveys \cite{Kobayashi-survey} and~\cite{Labourie-survey} for
references to some of the previous work on the existence of compact
Clifford-Klein forms.
 \end{defn}

We determine exactly which homogeneous spaces of $\SO(2,n)$ have a compact
Clifford-Klein form in the case where $n$~is even \see{CK-even}, and we have
almost complete results in the case where $n$~is odd \see{SU1m->complete}.
(We only consider homogeneous spaces $G/H$ in which $H$ is connected.) The
work leads to new examples of homogeneous spaces that have compact
Clifford-Klein forms, if $n$ is even \see{compact}. We also show that only
the obvious homogeneous spaces of $\SL(3,\real)$ have compact Clifford-Klein
forms \see{no-cpt}, and we study noncompact Clifford-Klein forms of finite
volume (see~\S\ref{finvol-sect}). 

\begin{notation} \label{SO2n-defn}
 We realize $\SO(2,n)$ as isometries of the indefinite form $\langle v \mid
v \rangle = v_1 v_{n+2} + v_2 v_{n+1} + \sum_{i=3}^n v_i^2$ on~$\real^{n+2}$
(for $v = (v_1,v_2,\ldots,v_{n+2}) \in \real^{n+2}$).
 Let $A$ be the subgroup consisting of the diagonal matrices in $\SO(2,n)$
whose diagonal entries are all positive, and let $N$ be the subgroup
consisting of the upper-triangular matrices in $\SO(2,n)$ with only $1$'s on
the diagonal.
 Thus, the Lie algebra of $AN$ is
 \begin{equation}
 \label{SO2n-AN}
 \Lie a + \Lie n =
 \bigset{
  \begin{pmatrix}
 t_1 & \phi & x  & \eta &0\\
 & t_2 & y &0 & -\eta \\
  &   & 0 & -y^T & -x^T \\
 && &-t_2& -\phi \\
 &&& &-t_1 \\
 \end{pmatrix}
 }{ {t_1,t_2, \phi,\eta \in \real, \atop x,y \in \real^{n-2}}}
 .
 \end{equation}
 Note that  the first two rows of any element of $\Lie a + \Lie n$ are
sufficient to determine the entire matrix.
 \end{notation}

Let us recall a construction of compact Clifford-Klein forms found by
Kulkarni \cite[Thm.~6.1]{Kulkarni} (see also
\cite[Prop.~4.9]{Kobayashi-properaction}). The subgroup $\SU(1,m)$,
embedded into $\SO(2,2m)$ in a standard way, acts properly and
transitively on the homogeneous  space $\SO(2,2m)/ \SO(1,2m)$. Therefore,
any co-compact lattice~$\Gamma$ in $\SO(1,2m)$ acts properly on
$\SO(2,2m)/\SU(1,m)$, and the quotient $\Gamma \backslash \SO(2,2m)/
\SU(1,m)$ is compact. Now let $H_{\SU} = \SU(1,m) \cap (AN)$. Then the
Clifford-Klein form $\Gamma \backslash \SO(2,2m)/ H_{\SU}$ is also compact,
since $\SU(1,m) / H_{\SU}$ is compact. (Similarly, Kulkarni also constructed
compact Clifford-Klein forms $\Lambda \backslash \SO(2,2m) / \SO(1,2m)$, by
letting $\Lambda$ be a co-compact lattice in $\SU(1,m)$.)

The following theorem demonstrates how to construct new examples of
compact Clifford-Klein forms $\Gamma \backslash \SO(2,2m) / H_B$. The
subgroup~$H_B$ of $\SO(2,2m)$ is obtained by deforming $H_{\SU}$, but
$H_B$ is almost never contained in any conjugate of $\SU(1,m)$.

\begin{thm} \label{compact}
 Assume that $G = \SO(2,2m)$.
 Let $B \colon \real^{2m-2} \to \real^{2m-2}$ be a linear transformation that
has no real eigenvalue. Set
 \begin{equation} \label{HB-defn}
 \Lie h_B=\bigset{
 \begin{pmatrix}
 t & 0 & x & \eta          & 0  \\
   &t  & B(x)    & 0             & -\eta \\
   &   & \dots
 \end{pmatrix}
 }{
 {x\in \real ^{2m-2}, 
 \atop
 t,\eta \in \real}
 } ,
 \end{equation}
 let $H_B$ be the corresponding closed, connected subgroup of~$G$,
 and let  $\Gamma$ be a co-compact lattice in $\SO(1,2m)$.
 Then
 \begin{enumerate}
 \item \label{HB-proper}
 the subgroup $\Gamma$ acts properly on $\SO(2,2m)/ H_B$;
 \item \label{HB-compact}
 the quotient $\Gamma \backslash \SO(2,2m)/ H_B$ is compact; and
 \item \label{HB-conj}
 $H_B$ is conjugate via $\On(2,2m)$ to a subgroup of $\SU(1,m)$ if and only
if for some $a, b \in \real$ {\upshape(}with $b \neq 0${\upshape)}, the
matrix of~$B$ with respect to some orthonormal basis of~$\real^{2m-2}$ is a
block diagonal matrix each of whose blocks is
 $\begin{pmatrix} a& b\\-b & a \end{pmatrix}$.
 \end{enumerate}
 Furthermore, by varying~$B$, one can obtain uncountably many pairwise
nonconjugate subgroups.
 \end{thm}

We recall that T.~Kobayashi \cite[Thm.~B]{Kobayashi-deformation} showed
that a co-compact lattice in $\SU(1, m)$ can be deformed to a discrete
subgroup~$\Lambda$, such that $\Lambda$ acts properly on $\SO(2,2m) /
\SO(1,2m)$ and  the quotient space  $\Lambda \backslash \SO(2,2m) /
\SO(1,2m)$ is compact, but $\Lambda$ is not contained in any conjugate of
$\SU(1,m)$. (This example is part of an extension of work of
W.~Goldman~\cite{Goldman}.) Note that Kobayashi created new compact
Clifford-Klein forms by deforming the discrete group while keeping the
homogeneous space $\SO(2,2m) / \SO(1,2m)$ fixed. In contrast, we deform
the homogeneous space $\SO(2,2m) /H_{\SU}$ to another homogeneous space
$\SO(2,2m) /H_B$ while keeping the discrete group~$\Gamma$ in $\SO(1,2m)$
fixed.

For even~$n$, we show that the Kulkarni examples and our deformations are
essentially the only interesting homogeneous spaces of $\SO(2,n)$ that have
compact Clifford-Klein forms. We assume that $H \subset AN$, because the
general case reduces to this \see{HcanbeAN}.

\begin{thm}[{cf.~Thm.~\ref{noCK}}] \label{CK-even}
 Assume that $G = \SO(2,2m)$. Let $H$ be a closed, connected subgroup
of $AN$, such that neither $H$ nor $G/H$ is compact.
 The homogeneous space $G/H$ has a compact Clifford-Klein form if and only
if either
 \begin{enumerate}
 \item \label{CK-even-SO1n}
 $H$ is conjugate to a co-compact subgroup of $\SO(1,2m)$; or
 \item \label{CK-even-SU1n}
 $H$ is conjugate to $H_B$, for some~$B$, as
described in Theorem~\ref{compact}.
 \end{enumerate}
 \end{thm}

It is conjectured \cite[1.4]{Kobayashi-deformation} that if $H$ is
reductive and $G/H$ has a compact Clifford-Klein form, then there exists
a reductive subgroup~$L$ of~$G$, such that $L$ acts properly on $G/H$,
and the double-coset space $L \backslash G/H$ is compact. Because there
is no such subgroup~$L$ in the case where $G = \SO(2,2m+1)$ and $H =
\SU(1,m)$ \see{SOodd/SU}, the following is a special case of the general
conjecture.

\begin{conj} \label{SU1m-conj}
 For $m \ge 1$, the homogeneous space  $\SO(2,2m+1)/\SU(1,m)$ does
not have a compact Clifford-Klein form.
 \end{conj}

If this conjecture is true, then, for odd~$n$, there is no interesting
example of a homogeneous space of $\SO(2,n)$ that has a compact
Clifford-Klein form.

 \begin{thm} \label{SU1m->complete}
  Assume that $G = \SO(2,2m+1)$, and that $G/\SU(1,m)$ does not have a
compact Clifford-Klein form. If $H$ is a closed, connected subgroup of~$G$,
such that neither $H$ nor $G/H$ is compact, then $G/H$ does not have a
compact Clifford-Klein form.
 \end{thm}

The main results of~\cite{OhWitte-CDS} list the homogeneous spaces of
$\SO(2,n)$ that admit a proper action of a noncompact subgroup of $\SO(2,n)$
(see~\S\ref{prelims}). Our proofs of Theorems~\ref{CK-even}
and~\ref{SU1m->complete} consist of case-by-case analysis to decide whether
each of these homogeneous spaces has a compact Clifford-Klein form. The
following proposition does not require such a detailed analysis, but is
obtained easily by combining theorems of Y.~Benoist \see{Ben-B} and
G.~A.~Margulis \see{MargulisTempered}.

\begin{prop} \label{no-cpt}
 Let $H$ be a closed, connected subgroup of $G = \SL(3,\real)$. If $G/H$ has
a compact Clifford-Klein form, then $H$ is either compact or co-compact.
 \end{prop}

The paper is organized as follows. Section~\ref{prelims} recalls some
definitions and results, mostly from \cite{OhWitte-CDS}.
Section~\ref{general} presents some general results on Clifford-Klein forms.
Section~\ref{examples} proves Theorem~\ref{compact}, the
new examples of compact Clifford-Klein forms. Section~\ref{SO2n-sect} proves
Theorems~\ref{CK-even} and~\ref{SU1m->complete}, the classification of
compact Clifford-Klein forms. Section~\ref{finvol-sect} discusses noncompact
Clifford-Klein forms of finite volume. Section~\ref{SL3-sect} briefly
discusses Clifford-Klein forms of homogeneous spaces of $\SL(3,\real)$,
proving Proposition~\ref{no-cpt}. An appendix presents a short proof of a
theorem of Benoist.

\begin{ack} This research was partially supported by grants from the
National Science Foundation (DMS-9623256 and DMS-9801136). 
 Many of the results of this paper were finalized during a visit of the
authors to the University of Bielefeld. We would like to  thank the
German-Israeli Foundation for Research and Development for financial
support that made the visit possible, and the mathematics faculty of the
University of Bielefeld for their hospitality that helped make the visit
so productive. D.~Witte would like to thank
the Tata Institute for Fundamental Research for providing a congenial
environment to carry out final revisions on the manuscript. Both authors are
thankful to Robert J.~Zimmer for his interest in this problem and his
encouragement.
 \end{ack}

\section{Cartan projections of subgroups of $\SO(2,n)$} \label{prelims}

Our study of the compact Clifford-Klein forms of homogeneous spaces of
$\SO(2,n)$ is based on a case-by-case analysis of the subgroups of
$\SO(2,n)$ that are not Cartan-decomposition subgroups (see
Defn.~\ref{CDS-defn}). We need to know not only what the subgroups are, but
also the image of each subgroup under the Cartan projection (see
Defn.~\ref{Cartproj-defn}). This information is provided
by~\cite{OhWitte-CDS}. In this section, we recall the relevant results,
notation, and definitions. (However, we have omitted a few of the more
well-known definitions that appear in \cite{OhWitte-CDS}.)

\begin{defn}[{\cite[Defn.~1.2]{OhWitte-CDS}}] \label{CDS-defn}
 Let $H$ be a closed, connected subgroup of~$G$. We say that $H$ is a
\emph{Cartan-decomposition subgroup} of~$G$ if there is a compact
subset~$C$ of~$G$, such that $CHC = G$. (Note that $C$ is only assumed to
be a sub\emph{set} of~$G$; it need not be a sub\emph{group}.)
 \end{defn}

\begin{notation}
 Fix an Iwasawa decomposition  $G = KAN$ and a corresponding Cartan
decomposition $G = K A^+ K$, where $A^+$ is the (closed) positive Weyl
chamber of~$A$ in which the roots occurring in the Lie algebra of~$N$ are
positive.  Thus, $K$ is a maximal compact subgroup, $A$ is the identity
component of a maximal split torus, and $N$ is a maximal unipotent subgroup.
 \end{notation}

\begin{defn}[Cartan projection] \label{Cartproj-defn}
 For each element~$g$ of~$G$, the Cartan decomposition $G = K A^+ K$
implies that there is an element~$a$ of~$A^+$ with $g \in K a K$. In fact,
the element~$a$ is unique, so there is a well-defined function $\mu
\colon G \to A^+$ given by $g \in K \, \mu(g) \, K$. The function $\mu$ is
continuous and proper (that is, the inverse image of any compact set
is compact). Some properties of the Cartan projection are discussed
in~\cite{Benoist}, \cite{Kobayashi-survey}, and~\cite{OhWitte-CDS}.
 \end{defn}

\begin{notation} 
 For subsets $X,Y \subset A^+$, we write $X \approx Y$ if there is a
compact subset $C$ of~$A$ with $X \subset YC$ and $Y \subset XC$.
 \end{notation}

\begin{notation} \label{SO2n-rho}
 Define a representation 
 $$\rho \colon \SO(2,n) \to \SL \bigl( \real^{n+2} \wedge \real^{n+2} \bigr)
 \mbox{ by }
 \rho(g) = g \wedge g .$$
 For functions $f_1,f_2 \colon \real^+ \to \real^+$, and a subgroup~$H$
of~$\SO(2,n)$, we write 
 $\mu(H) \approx \muH{f_1(\|h\|)}{f_2(\|h\|)}$ if, for every sufficiently
large $C > 1$, we have
 $$ \mu(H) \approx \bigset{a \in A^+}
 { C^{-1} f_1 \bigl( \|a\|
\bigr) \le \|\rho(a)\| \le C f_2 \bigl( \|a\| \bigr)} ,$$
 where $\|h\|$ denotes the norm of the linear transformation~$h$. (If
$f_1$~and~$f_2$ are monomials, or other very tame functions, then it does not
matter which particular norm is used, because all norms are equivalent up to
a bounded factor.)
 In particular, $H$ is a Cartan-decomposition subgroup of~$\SO(2,n)$ if and
only if $\mu(H) \approx \muH{\|h\|}{\|h\|^{2}}$ (see
\cite[Lem.~5.7]{OhWitte-CDS}).
 \end{notation}

\begin{rem} \label{H'inCHC}
 Let $H$ and~$H'$ be closed, connected subgroups of $\SO(2,n)$. Suppose that
 $$f_1,f_2,f_1',f_2' \colon \real^+ \to \real^+$$
 satisfy
 $$ \mu(H) \approx \muH{f_1(\|h\|)}{f_2(\|h\|)}
 \text{\qquad and\qquad} 
 \mu(H') \approx \muH{f_1'(\|h\|)}{f_2'(\|h\|)} .$$
 If, for all large $t \in \real^+$, we have $f_1(t) \le f_1'(t) \le f_2'(t)
\le f_2(t)$, then there is a compact subset~$C$ of~$G$, such that $H'
\subset CHC$.
 \end{rem}

\begin{notation}[cf.~\ref{SO2n-AN}]
 Assume that $G = \SO(2,n)$. For every $h \in \Lie n$, there exist unique
$\phi_h,\eta_h \in \real$ and $x_h,y_h \in \real^{n-2}$, such that
 $$ h = \begin{pmatrix}
 0 & \phi_h & x_h & \eta_h & 0 \\
   &  0     & y_h &  0     & - \eta_h \\
   &        & \dots \\
 \end{pmatrix}
 .$$
 \end{notation}

\begin{notation}
 We let $\alpha$ and~$\beta$ be the simple real roots of $\SO(2,n)$,
defined by $\alpha(a) = a_1/a_2$ and $\beta(a) = a_2$, for an element~$a$
of~$A$ of the form
 $$ a= \diag (a_1, a_2 , 1,1,\ldots,1,1 , a_2^{-1} ,
a_1^{-1}) . $$
 Thus, 
 \begin{itemize}
 \item the root space $\Lie u_\alpha$ is the $\phi$-axis in~$\Lie n$, 
 \item the root space $\Lie u_\beta$ is the $y$-subspace in~$\Lie n$, 
 \item the root space $\Lie u_{\alpha+\beta}$ is the $x$-subspace in~$\Lie n$, and
 \item  the root space $\Lie u_{\alpha+2\beta}$ is the $\eta$-axis
in~$\Lie n$.
 \end{itemize}
 \end{notation}

We now reproduce a string of results from~\cite{OhWitte-CDS} that describe
the subgroups of $\SO(2,n)$ that are not Cartan-decomposition subgroups, and
also describe the image of each subgroup under the Cartan projection.

\begin{thm}[{\cite[Thm.~5.5 and Prop.~5.8]{OhWitte-CDS}}] \label{HinN}
 Assume that $G = \SO(2,n)$. A closed, connected subgroup~$H$ of~$N$ is not a
Cartan-decomposition subgroup of~$G$ if and only if either
 \begin{enumerate}

\item \label{HinN-dim1}
 $\dim H \le 1$, in which case $\mu(H) \approx \muH{\|h\|^p}{\|h\|^p}$, for
some $p \in \{0,1,3/2,2\}$; or

 \item \label{HinN-<xy>not1}
 for every nonzero element~$h$ of~$\Lie h$, we have $\phi_h = 0$ and $\dim
\langle x_h,y_h \rangle \neq 1$, in which case $\mu(H)
\approx \muH{\|h\|^2}{\|h\|^2}$; or

 \item \label{HinN-<xy>=1}
 for every nonzero element~$h$ of~$\Lie h$, we have $\phi_h = 0$ and $\dim
\langle x_h,y_h \rangle = 1$, in which case $\mu(H) \approx
\muH{\|h\|}{\|h\|}$; or

 \item \label{HinN-b2-2a}
 there exists a subspace~$X_0$ of~$\real^{n-2}$, $b \in X_0$, $c \in
X_0^\perp$, and $p \in \real$ with $\|b\|^2 - \|c\|^2 - 2p < 0$, such that
for every element of~$\Lie h$, we have $y_h= 0$, $x_h \in \phi_h c + X_0$, and
$\eta_h = p\phi_h + b \cdot x_h$ {\upshape(}where $b \cdot x_h$ denotes the
Euclidean dot product of the vectors $b$ and~$x_h$
in~$\real^{n-2}${\upshape)}, in which case $\mu(H) \approx
\muH{\|h\|}{\|h\|}$.
 \end{enumerate}
 \end{thm}

\begin{thm}[{\cite[Cor.~6.2]{OhWitte-CDS}}] \label{SO2n-semi-notCDS}
 Assume that $G = \SO(2,n)$. Let $H$ be a closed, connected subgroup
of~$AN$, such that $H = (H \cap A) \semiprod (H \cap N)$,
and $H \not\subset N$. The subgroup~$H$ is not a Cartan-decomposition
subgroup of~$G$ if and only if either
 \begin{enumerate}

 \item \label{SO2n-semi-notCDS-1D}
 $H = H \cap A$ is a one-dimensional subgroup of~$A$; or

 \item \label{SO2n-semi-notCDS-<xy>not1}
 $H \cap A = \ker\alpha$, and we have $\phi_h = 0$ and
$\dim \langle x_h,y_h \rangle \neq 1$ for every nonzero element~$h$ of~$\Lie
h \cap \Lie n$, in which case  $\mu(H) \approx \muH{\|h\|^2}{\|h\|^2}$;
or

 \item \label{SO2n-semi-notCDS-y=0}
 $H \cap A = \ker\beta$, and we have $\phi_h =0$, $y_h = 0$, and $x_h \neq
0$ for every nonzero element~$h$ of $\Lie h \cap \Lie n$, in which case
$\mu(H) \approx \muH{\|h\|}{\|h\|}$; or

 \item \label{SO2n-semi-notCDS-x=0}
 $H \cap A = \ker(\alpha+\beta)$, and we have $\phi_h=0$, $x_h = 0$,
and $y_h \neq 0$, for every nonzero element~$h$ of $\Lie h \cap \Lie n$, in
which case $\mu(H) \approx \muH{\|h\|}{\|h\|}$; or

 \item \label{SO2n-semi-notCDS-b2-2p}
 $H \cap A = \ker \beta$, and there exist a
subspace~$X_0$ of~$\real^{n-2}$, $b \in X_0$, $c \in X_0^\perp$, and $p
\in \real$, such that $\|b\|^2 - \|c\|^2 - 2p < 0$, and we have $y_h = 0$,
$x_h \in \phi_h c + X_0$, and $\eta_h = p\phi + b \cdot x_h$ for every $h \in
\Lie h \cap \Lie n$, in which case $\mu(H) \approx \muH{\|h\|}{\|h\|}$; or

 \item \label{SO2n-semi-notCDS-phiy}
 $H \cap A = \ker(\alpha-\beta)$, $\dim H = 2$,  and there are
${\hat \phi} \in \Lie u_\alpha$ and ${\hat y} \in \Lie
u_{\beta }$, such that ${\hat \phi} \neq 0$, ${\hat y} \neq 0$, and $\Lie
h \cap \Lie n = \real({\hat \phi} + {\hat y})$, in which case
 $\mu(H) \approx \muH{\|h\|^{3/2}}{\|h\|^{3/2}}$; or

 \item \label{SO2n-semi-notCDS-dim2}
 $H \cap A =\ker\beta$, $\dim H = 2$, and we have $y_h = 0$
and $\|x_h\|^2 \neq -2 \phi_h \eta_h$ for every $h \in H$, in which case
$\mu(H) \approx \muH{\|h\|}{\|h\|}$; or

 \item \label{SO2n-semi-notCDS-root}
 there is a positive root~$\omega$ and a one-dimensional
subspace~$\Lie t$ of~$\Lie a$, such that $\Lie h \cap \Lie n \subset \Lie
u_\omega$, $\Lie h = \Lie t + (\Lie h \cap \Lie n)$, and
Proposition~\ref{rootsemi-CDS} implies that $H$ is not a
Cartan-decomposition subgroup.

 \end{enumerate}
 \end{thm}

Not every connected subgroup of $AN$ is conjugate to a subgroup of the form
$T \semiprod U$, where $T \subset A$ and $U \subset N$. The following
definition and lemma describe how close we can come to this ideal situation.

\begin{defn} \label{compatible}
 Let us say that a subgroup~$H$ of $AN$ is \emph{compatible} with~$A$ if
$H \subset T U C_N(T)$, where $T = A \cap (HN)$, $U = H \cap N$, and
$C_N(T)$ denotes the centralizer of~$T$ in~$N$. 
 \end{defn}

\begin{lem}[{\cite[Lem.~2.3]{OhWitte-CDS}}] \label{conj-to-compatible}
 If $H$ is a closed, connected subgroup of~$AN$, then $H$ is conjugate,
via an element of~$N$, to a subgroup that is compatible with~$A$.
 \end{lem}

\begin{thm}[{\cite[Cor.~6.4]{OhWitte-CDS}}]
\label{SO2n-notsemi-notCDS}
 Assume that $G = \SO(2,n)$. Let $H$ be a closed, connected subgroup of~$AN$
that is compatible with~$A$ \see{compatible}, and assume that $H \neq (H
\cap A) \semiprod (H \cap N)$. Then there is a positive root~$\omega$, and a
one-dimensional subspace~$\Lie x$ of $(\ker \omega) + {\Lie u}_\omega$, such
that $\Lie h = \Lie x + (\Lie h \cap \Lie n)$.

 If $H$ is not a Cartan-decomposition subgroup of~$G$, then either:
 \begin{enumerate}
 \item \label{SO2n-notsemi-notCDS-a,a+b}
  $\omega =\alpha$ and $\Lie h \cap \Lie n \subset \Lie u_{\alpha+\beta}$,
in which case
 $\mu(H) \approx \muH{\|h\|}{\|h\|^2/(\log \|h\|)}$; or
 \item \label{SO2n-notsemi-notCDS-a,a+2b}
  $\omega =\alpha$ and $\Lie h \cap \Lie n \subset \Lie u_{\alpha+2\beta}$,
in which case
 $\mu(H) \approx \muH{\|h\|^2/(\log \|h\|)^2}{\|h\|^2}$; or
 \item \label{SO2n-notsemi-notCDS-a+2b,a}
  $\omega ={\alpha + 2\beta}$ and $\Lie h \cap \Lie n \subset \Lie
u_{\alpha}$, in which case
 $\mu(H) \approx \muH{\|h\|^2/(\log \|h\|)^2}{\|h\|^2}$; or
 \item \label{SO2n-notsemi-notCDS-a+2b,b}
  $\omega ={\alpha + 2\beta}$ and either $\Lie h \cap \Lie n \subset \Lie
u_{\beta}$ or~$\Lie h \cap \Lie n \subset \Lie u_{\alpha+\beta}$, in which
case
 $\mu(H) \approx \muH{\|h\|}{\|h\|^2/(\log \|h\|)}$; or
 \item \label{SO2n-notsemi-notCDS-b(a+b),w(a+2b)}
  $\omega \in\{\beta, \alpha + \beta\}$, $\Lie h \cap \Lie n \subset \Lie
u_{\omega} + \Lie u_{\alpha+2\beta}$, and $\Lie h \cap \Lie u_\omega = \Lie
h \cap \Lie u_{\alpha+2\beta} = 0$, in which case
 $\mu(H) \approx \muH{\|h\|}{\|h\|^{3/2}}$; or 
 \item \label{SO2n-notsemi-notCDS-b(a+b),c(a+2b)}
  there is a root~$\gamma$ with $\{\omega, \gamma\} = \{\beta, \alpha +
\beta\}$, $\Lie h \cap \Lie n \subset \Lie u_{\gamma} + \Lie
u_{\alpha+2\beta}$, and $\Lie h \cap \Lie u_{\alpha+2\beta} = 0$, in which
case
 $\mu(H) \approx \muH{\|h\|}{\|h\| (\log\|h\|)^2}$; or 
 \item \label{SO2n-notsemi-notCDS-b(a+b),a+2b}
  $\omega \in \{\beta, \alpha + \beta\}$ and $\Lie h \cap \Lie n = \Lie
u_{\alpha+2\beta}$, in which case
 $\mu(H) \approx \muH{\|h\| (\log\|h\|)}{\|h\|^2}$; or
 \item \label{SO2n-notsemi-notCDS-a+b,a}
   $\omega = \alpha + \beta$ and $\Lie h \cap \Lie n = \Lie u_\alpha$, in
which case
 $\mu(H) \approx \muH{\|h\| (\log\|h\|)}{\|h\|^2}$.
 \end{enumerate}
 \end{thm}

\begin{prop}[{cf.~\cite[Prop.~3.17, Cor.~3.18]{OhWitte-CDS}}]
\label{rootsemi-CDS}
 Assume that $G = \SO(2,n)$. Let $H$ be a closed, connected subgroup
of~$AN$, such that there is a positive root~$\omega$ and a one-dimensional
subspace~$\Lie t$ of~$\Lie a$, such that $0 \neq \Lie h \cap \Lie n \subset
\Lie u_\omega$, and $\Lie h = \Lie t + (\Lie h \cap \Lie n)$.

Let $\gamma$ be the positive root that is perpendicular to~$\omega$, so
$\{\omega,\gamma\}$ is either $\{\alpha, \alpha+2\beta\}$ or $\{\beta,
\alpha+\beta\}$. If $\Lie t = \ker \omega$, then $H$ is a
Cartan-decomposition subgroup of~$G$. Otherwise, there is a real number~$p$,
such that $\Lie t = \ker( p \omega + \gamma)$.
 \begin{itemize}
 \item If $|p| \ge 1$, then $H$ is a Cartan-decomposition subgroup of~$G$.
 \item If $|p| < 1$, and $\omega \in \{\alpha, \alpha+2\beta\}$, then
$\mu(H) \approx \muH{\|h\|^{2/(1+|p|)}}{\|h\|^2}$.
 \item If $|p| < 1$, and $\omega \in \{\beta, \alpha+\beta\}$, then
$\mu(H) \approx \muH{\|h\|}{\|h\|^{1+|p|}}$.
 \end{itemize}
 \end{prop}

For ease of reference, we now collect a few miscellaneous facts.

\begin{lem}[{\cite[Lem.~2.8]{OhWitte-CDS}}] \label{HN=AN}
 Let $H$ be a closed, connected subgroup of~$AN$.
 If $\dim H - \dim(H \cap N) \ge \Rrank G$, then $H$ contains a
conjugate of~$A$, so $H$ is a Cartan-decomposition subgroup.
 \end{lem}

\begin{rem} \label{mu(SO1n)}
 We realize $\SO(1,n)$ as the stabilizer of the vector
 $$(0,1,0,0,\ldots,0,-1,0) .$$
 Thus, the Lie algebra of $\SO(1,n) \cap AN$ is
 $$\bigset{
 \begin{pmatrix}
 t & \phi & x & \phi & 0 \\
   &  0   & 0 &   0  & - \phi \\
   &      & \dots  \\
 \end{pmatrix}
 }
 { {t, \phi \in \real \atop x \in \real^{n-2}}
 }
 ,$$
 that is, it is of type~\fullref{SO2n-semi-notCDS}{-b2-2p}, with $X_0 =
\real^{n-2}$, $b = 0$, $c = 0$, and $p = 1$. Therefore, we see that $\mu
\bigl( \SO(1,n) \bigr) \approx \muH{\|h\|}{\|h\|}$.
 \end{rem}

\begin{prop}[{cf.~\cite[Case~3 of the pf.\ of Thm.~6.1]{OhWitte-CDS}}]
\label{SO1n-whenconj}
 Assume that $G = \SO(2,n)$. Suppose that $H$ is a closed, connected
subgroup of~$AN$, such that $H = (H \cap A) \semiprod (H \cap N)$. Assume
that there exist a subspace~$X_0$ of~$\real^{n-2}$, vectors $b \in X_0$ and
$c \in X_0^\perp$, and a real number~$p$, such that, for every $h \in H \cap
N$, we have
 $y_h = 0$, $x_h \in \phi_h c + X_0$, and $\eta_h = p\phi_h + b \cdot x_h$.
 If $\|b\|^2 - \|c\|^2 - 2p < 0$, then $H$ is conjugate to a subgroup of
$\SO(1,n)$.
 \end{prop}

\section{General results on compact Clifford-Klein forms} \label{general}

Before beginning our study of the specific group $\SO(2,n)$, let us state
some general results on compact Clifford-Klein forms. Recall that $G$ is a
Zariski-connected, almost simple, linear Lie group.

The Calabi-Markus phenomenon asserts that if $H$ is a Cartan-decomposition
subgroup of~$G$, then no closed, noncompact subgroup of~$G$ acts properly on
$G/H$ (cf.~\cite[pf.~of Thm.~A.1.2]{Kulkarni}). The following well-known
fact is a direct consequence of this observation. 

\begin{lem} \label{CDS->notess}
 Let $H$ be a Cartan-decomposition subgroup of~$G$. 
 Then $G/H$ does not have a compact Clifford-Klein form, unless $G/H$ itself
is compact.
 \end{lem}

By combining this lemma with the following proposition, we see that if $G/H$
has a compact Clifford-Klein form, and $\Rrank G = 1$, then either $H$ or
$G/H$ is compact.

\begin{prop}[{\cite[Prop.~1.8]{OhWitte-CDS},
\cite[Lem.~3.2]{Kobayashi-isotropy}}] \label{Rrank1-CDS}
 Assume that $\Rrank G = 1$. A closed, connected subgroup~$H$ of~$G$ is a
Cartan-decomposition subgroup if and only if $H$ is noncompact.
 \end{prop}

Lemma~\ref{CDS->notess} can be obtained as a special case of
Theorem~\fullref{noncpct-dim}{-d(L)>d(H)} by letting $L = G$. The
generalization is very useful.

\begin{notation}[{cf.~\cite[(2.5), \S5]{Kobayashi-properaction}}]
\label{d(H)-defn}
 For any connected Lie group~$H$, let $d(H) = \dim H - \dim K_H$, where
$K_H$ is a maximal compact subgroup of~$H$. This is well defined, because
all the maximal compact subgroups of~$H$ are conjugate \cite[Thm.~XV.3.1,
p.~180--181]{Hochschild-Lie}. (This concept originated with K.~Iwasawa
\cite[p.~533]{Iwasawa}, who called $d(H)$ the ``characteristic index"
of~$H$.) Note that if $H \subset AN$, then $d(H) = \dim H$, because $AN$
has no nontrivial compact subgroups.
 \end{notation}

\begin{thm}[{Kobayashi, cf.~\cite[Cor.~5.5]{Kobayashi-properaction} and
\cite[Thm.~1.5]{Kobayashi-necessary}}] \label{noncpct-dim}
 Let $H$ and~$L$ be closed, connected subgroups of~$G$, and assume that
there is a compact subset~$C$ of~$G$, such that $L \subset CHC$.
 \begin{enumerate}
 \item \label{noncpct-dim-d(L)>d(H)}
 If $d(L) > d(H)$, then $G/H$ does not have a compact
Clifford-Klein form.
 \item \label{noncpct-dim-d(L)=d(H)}
 If $d(L) = d(H)$, and $G/H$ has a compact
Clifford-Klein form, then $G/L$ also has a compact Clifford-Klein form.
 \item \label{noncpct-dim-d(L')+d(H)=d(G)}
 If there is a closed subgroup~$L'$ of~$G$,
such that $L'$ acts properly on $G/H$, $d(H) + d(L') = d(G)$, and there
is a co-compact lattice~$\Gamma$ in~$L'$, then $G/H$ has a compact
Clifford-Klein form, namely, the quotient $\Gamma \backslash G/H$ is
compact.
 \end{enumerate}
 \end{thm}

\begin{proof}[Comments on the proof]
 Kobayashi assumed that $H$ is reductive, but the same proof works with
only minor changes. From Lemma~\ref{HcanbeAN}, we see that we may
assume that $H \subset AN$. Then, from the Iwasawa decomposition $G = KAN$,
it is immediate that the homogeneous space $G/H$ is homeomorphic to the
cartesian product
 $K \times (AN/H)$.
 Because $AN$ is a simply connected, solvable Lie group and $H$ is a
connected subgroup, the homogeneous space $AN/H$ is homeomorphic to a
Euclidean space~$\real^M$ (cf.~\cite[Prop.~11.2]{MostowFSS}). Therefore,
we see that $G/H$ has the same homotopy type as $K$ ($= K/(H \cap K)$,
because $H \cap K$ is trivial). Thus, the proof of
\cite[Lem.~5.3]{Kobayashi-properaction} goes through essentially
unchanged in the general setting. This yields a general version of
\cite[Cor.~5.5]{Kobayashi-properaction}, from which general versions of
\cite[Thm.~1.5]{Kobayashi-necessary} and
\cite[Thm.~4.7]{Kobayashi-properaction} follow. Our
conclusion~\pref{noncpct-dim-d(L)>d(H)} is the natural generalization of
\cite[Thm.~1.5]{Kobayashi-necessary}.
Conclusion~\pref{noncpct-dim-d(L')+d(H)=d(G)} is the natural
generalization of \cite[Thm.~4.7]{Kobayashi-properaction}, and
conclusion~\pref{noncpct-dim-d(L)=d(H)} is proved similarly.
 \end{proof}

The following useful theorem of G.~A.~Margulis is used in the proof of
Proposition~\ref{tempered}. 

\begin{defn}[{cf.~\cite[Defn.~2.2, Rmk.~2.2]{Margulis}}] \label{tempered-def}
 A closed subgroup~$H$ of~$G$ is \emph{$(G,K)$-tempered} if there exists a
(positive) function $q \in L^1(H)$, such that, for every non-trivial,
irreducible, unitary representation~$\pi$ of~$G$ with a $K$-fixed unit
vector~$v$, we have $|\langle \pi (h) v , v \rangle  | \leq q(h)$ for all $h
\in H$. 

(We remark that, because $\pi$ is irreducible, the $K$-invariant vector~$v$
is unique, up to a scalar multiple \cite[Thm.~8.1]{KnappExamples}.)
 \end{defn}

\begin{thm}[{Margulis \cite[Thm.~3.1]{Margulis}}] \label{MargulisTempered}
 If $H$ is a closed, noncompact, $(G,K)$-tempered subgroup of~$G$, then
$G/H$ does not have a compact Clifford-Klein form.
 \end{thm}

\begin{prop} \label{tempered}
 If $H$ is a closed, noncompact one-parameter subgroup of a connected,
simple, linear Lie group~$G$, then $G/H$ does not have a compact
Clifford-Klein form.
 \end{prop}

\begin{proof}
 Suppose that $G/H$ does have a compact Clifford-Klein form.

 Assume for the moment that $\Rrank G = 1$. Then $H$ is a
Cartan-decomposition subgroup of~$G$ \see{Rrank1-CDS}, so we see from
Lemma~\ref{CDS->notess} that $G/H$ must be compact. But the dimension of
every connected, co-compact subgroup of~$G$ is at least $d(G)$ \cite[(1.2),
p.~263]{GotoWang}, and $d(G) \ge 2$. This contradicts the fact that $\dim H
= 1$. Thus, we now know that $\Rrank G \ge 2$.

\setcounter{case}{0}

 \begin{case}
 Assume that $H$ is unipotent.
 \end{case}
 The Jacobson-Morosov Lemma \cite[Thm.~17(1), p.~100]{Jacobson}
implies that there exists a connected, closed subgroup $L$ of~$G$ that is
locally isomorphic to $\SL(2,\real)$ (and has finite center), and
contains~$H$. Then $H$ is a Cartan-decomposition subgroup of~$L$
\see{Rrank1-CDS}, and $d(L) = 2 > 1 = d(H)$, so
Theorem~\fullref{noncpct-dim}{-d(L)>d(H)} applies.

 \begin{case} Assume that $H$ is not unipotent. \end{case}
 We may assume that $H=\{\, a_t u_t\mid t\in \real \,\}$ where $a_t\in A$ is
a semisimple one parameter subgroup and $u_t\in N$ is a unipotent one
parameter subgroup such that $a_t$ commutes with $u_t$
(cf.~\ref{HcanbeAN} and~\ref{conj-to-compatible}). It
is well known (cf.~\cite[\S3, p.~140]{KatokSpatzier}, where a stronger
result is obtained by combining \cite[Cor.~7.2 and \S7]{Howe} with
\cite[Thm.~2.4.2]{Cowling}) that there are constants $C>0$ and $p>0$ such
that, for any non-trivial irreducible unitary representation~$\rho$ with
a $K$-fixed unit vector, say~$v$, we have
 $$| \langle \rho (g) v, v \rangle| \leq C \exp \bigl(-p d(e, g) \bigr)
 \mbox{ \qquad for $g \in G$} $$
where $d$ is some bi-$K$-invariant Riemannian metric on $G$. We may
assume that $d(e, a_t)=|t|$ with a suitable parameterization. Since the
growth of a unipotent one parameter subgroup is logarithmic while that of
semisimple one parameter subgroup is linear, we can find a large $T$ such
that
 $d(e, a_t u_t) \geq \frac{1}{2} d(e, a_t)= |t|/2 $ for all $t >T$.
 Hence
 $$| \langle \rho (a_t u_t) v, v \rangle| \leq C \exp \left(-\frac {p
|t|}2 \right)
 \mbox{ \qquad for $t > T$} .$$
 Since the function $ \exp (-p|t|/2)$ is in $L^1(\real)$, it follows that
$H$ is $(G,K)$-tempered. Therefore, Theorem~\ref{MargulisTempered}
implies that $G/H$ does not have compact Clifford-Klein forms.
 \end{proof}

We use the following well-known lemma to reduce the study of compact
Clifford-Klein forms of $G/H$ to the case where $H \subset AN$. We remark
that the proof of the lemma is constructive. For example, replace $H$ by a
conjugate, so that $\overline{H} \cap AN$ is co-compact in~$\overline{H}$,
where $\overline{H}$ is the Zariski closure of~$\overline{H}$, and choose a
maximal compact subgroup~$\overline{C}$ of~$\overline{H}^\circ$. Then write
$\overline{C} = C_1 C_2$, where $C_1$ is a maximal compact subgroup of~$H$,
and $C_2$ is contained in the Zariski closure of $\Rad H$. Finally, let $H'
= (HC_2) \cap (AN)$.

\begin{lem}[{cf.~\cite[Lem.~2.9]{OhWitte-CDS}}] \label{HcanbeAN}
 Let $H$ be a closed, connected subgroup of~$G$. Then there is a closed,
connected subgroup~$H'$ of~$G$ and compact, connected subgroups~$C_1$
and~$C_2$ of~$G$, such that
 \begin{enumerate}
 \item $H'$ is conjugate to a subgroup of~$AN$;
 \item $\dim H' = d(H)$ {\upshape(}see Notation~\ref{d(H)-defn}{\upshape)};
 \item $C_2 H = C_1 C_2 H'$;
 \item $C_1 \subset H$, $C_2$ is abelian, $C_1$ centralizes~$C_2$, and $C_2$
normalizes both~$H$ and~$H'$; and
 \item  if $H/\Rad H$ is compact, then $C_1$ normalizes~$H'$.
 \end{enumerate}
 Moreover, from~\fullref{noncpct-dim}{-d(L)=d(H)}, we know that the
homogeneous space $G/H$ has a compact Clifford-Klein form if and only if
$G/H'$ has a compact Clifford-Klein form.
 \end{lem}

We now recall a fundamental result of Benoist and Kobayashi.

\begin{thm}[{Benoist \cite[Prop.~1.5]{Benoist}, Kobayashi
\cite[Cor.~3.5]{Kobayashi-criterion}}] \label{proper}
 Let $H_1$ and~$H_2$ be closed subgroups of~$G$. The subgroup~$H_1$ acts
properly on~$G/H_2$ if and only if, for every compact subset~$C$ of~$A$,
the intersection $\bigl( \mu(H_1) C \bigr) \cap \mu(H_2)$ is compact.
 \end{thm}

\begin{lem} \label{Gconn}
 Let $G^\circ$ be the identity component of~$G$, let $H$ be a closed,
connected subgroup of~$G$, and let $\Gamma$ be a discrete subgroup of~$G$.
Then:
 \begin{enumerate}
 \item \label{Gconn-proper}
 $\Gamma$ acts properly on $G/H$ if and only if $\Gamma \cap G^\circ$
acts properly on $G^\circ/H$.
 \item \label{Gconn-cpct}
 $\Gamma \backslash G / H$ is compact if and only if  $(\Gamma \cap
G^\circ) \backslash G^\circ / H$ is compact.
 \end{enumerate}
 \end{lem}

\begin{proof}
 \pref{Gconn-proper} Because every element
of the Weyl group of~$G$ has a representative in~$G^\circ$
\cite[Cor.~14.6]{BorelTits}, we see that $G$ and~$G^\circ$ have the same
positive Weyl chamber~$A^+$, and the Cartan projection $G^\circ \to A^+$
is the restriction of the Cartan projection $G \to A^+$. Thus, the
desired conclusion is immediate from Corollary~\ref{proper}.

\pref{Gconn-cpct} This is an easy consequence of the fact that $G/G^\circ$ is
finite \cite[Appendix]{MostowSolv}.
 \end{proof}

\section{New examples of compact Clifford-Klein forms}
\label{examples}

\begin{proof}[{\bf Proof of Theorem~\ref{compact}}]
 From Remark~\ref{mu(SO1n)}, we have $\mu\bigl( \SO(1,2m) \bigr) \approx
\muH{\|h\|}{\|h\|}$. From Theorem~\fullref{SO2n-semi-notCDS}{-<xy>not1}, we
have $\mu(H_B) \approx \muH{\|h\|^2}{\|h\|^2}$. Thus,
Theorem~\ref{proper} implies that $\SO(1,2m)$ acts properly on
$\SO(2,2m)/H_B$. We have $d\bigl( \SO(2,2m) \bigr) = 4m$, $d \bigl( \SO(1,2m)
\bigr) = 2m$, and $\dim(H_B) = 2m$, so $d\bigl( \SO(2,2m) \bigr) = d \bigl(
\SO(1,2m) \bigr) + \dim(H_B)$. Therefore, conclusions \pref{HB-proper}
and~\pref{HB-compact} follow from
Theorem~\fullref{noncpct-dim}{-d(L')+d(H)=d(G)}.

To show conclusion~\pref{HB-conj}, suppose that $g\Lie h_B g^{-1} \subset
\Lie{su}(1,m)$ for some $g\in SO(2,2m)$. Because all maximal split tori in
$\SU(1,m)$ are conjugate, we may assume that $g$ normalizes $\ker\alpha$. In
fact, because $\SU(1,m)$ has an element (namely, the nontrivial element of
the Weyl group) that inverts $\ker \alpha$, we may assume that $g$
centralizes $\ker \alpha$. Thus, $g$ is a block diagonal matrix with $R$,
$S$ and $J(R^T)^{-1}J$ on the diagonal, where $R \in \GL(2, \real)$, $S
\in \On(2m-2)$, and $J=\begin{pmatrix} 0& 1\\ 1 & 0 \end{pmatrix}$.
Conjugating by $I_2 \times S \times I_2$ amounts to choosing a different
orthonormal basis for~$\real^{2m-2}$, so we may assume that $S$ is trivial.
Now, the assumption that $g\Lie h_B g^{-1} \subset \Lie{su}(1,m)$ implies
that
 $R\begin{pmatrix} x \\ B(x) \end{pmatrix}$
 is of the form
 $$\begin{pmatrix}
 z_1&-z_2 & z_3 &-z_4 & \hdots & z_{2m-3} &-z_{2m-2} \\
 z_2& z_1 & z_4 & z_3 & \hdots & z_{2m-2} & z_{2m-3}
 \end{pmatrix}
 .$$
  A direct calculation, setting 
 $$ R_{1,1} x_{2i-1} + R_{1,2} B(x)_{2i-1} = R_{2,1} x_{2i} + R_{2,2}
B(x)_{2i} $$
 and
 $$ R_{1,1} x_{2i} + R_{1,2} B(x)_{2i} = -\bigl( R_{2,1} x_{2i-1} + R_{2,2}
B(x)_{2i-1} \bigr) $$
 for all $i \in \{1,\ldots,m-1\}$, now establishes that $B$ is block
diagonal as described in conclusion~\pref{HB-conj}.

All that remains is to show that there are uncountably many nonconjugate
subgroups of the form~$H_B$. Let $L$ be the group of all block diagonal
matrices with $R$, $S$ and $J(R^T)^{-1}J$ on the diagonal, where $R \in
\GL(2, \real)$, $S \in \On(2m-2)$, and $J=\begin{pmatrix} 0& 1\\ 1 & 0
\end{pmatrix}$, and let $\mathcal{B}$ be the set of all matrices in
$\GL_{2m-2}(\real)$ that have no real eigenvalues, so $\mathcal{B}$ is an
open subset of $\GL_{2m-2}(\real)$. Because $B \mapsto H_B$ is injective,
the action of~$L$ by conjugation on $\{\, H_B \mid B \in \mathcal{B} \,\}$
yields an action of~$L$ on~$\mathcal{B}$. By arguing as in the proof
of~\pref{HB-conj}, we see that if $B_1$ and~$B_2$ are in different
$L$-orbits in~$\mathcal{B}$, then $H_{B_1}$ is not conjugate to~$H_{B_2}$.
Thus, it suffices to show that there are uncountably many $L$-orbits
on~$\mathcal{B}$.

 If $2m \ge 6$, then
 $$ \dim L = 4 + \frac{1}{2} (2m-2)(2m-3) < (2m-2)^2 = \dim \mathcal{B} ,$$
 so each $L$-orbit on~$\mathcal{B}$ has measure zero. Thus, obviously, there
are uncountably many $L$-orbits.

Now assume that $m = 2$. In this case, for each $B \in \mathcal{B}$, there
exists $B' \in \mathcal{B}$, such that
 $$ \bigset{ \begin{pmatrix} x \\ B(x) \end{pmatrix}} {x \in \real^2} = 
\bigset{ \begin{pmatrix} v & B'(v) \end{pmatrix}} {v \in \real^2} ,$$
 where $v$ and $B'(v)$ are considered as column vectors.
 Note that the centralizer of~$B'$ in $\GL(2,\real)$ contains a
2-dimensional connected subgroup.
 Thus, $B'$ is centralized by a nontrivial connected subgroup of
$\PGL(2,\real)$, so we see that the normalizer $N_L(H_B)$ contains a
2-dimensional subgroup (consisting of block diagonal matrices with $R$,
$\Id$ and $J(R^T)^{-1}J$ on the diagonal, for $R$ centralizing~$B'$), so
 $\dim L - \dim \bigl( N_L(H_B) \bigr) \le 3 < 4 = \dim \mathcal{B}$.
 Therefore, as in the previous case, each $L$-orbit on~$\mathcal{B}$ has
measure zero, so there are uncountably many $L$-orbits.
 \end{proof}

\begin{rem}
 The subgroups~$H_B$ of Theorem~\ref{compact} are not all isomorphic (unless
$m = 2$). For example,
 let $m = 3$ and let 
 $$ B = \begin{pmatrix}
  0&1&0&1 \\
 -1&0&1&0 \\
  0&1&0&0 \\
  1&0&0&0 \\
 \end{pmatrix}
 .
 $$
 The characteristic polynomial of~$B$ is $\det(\lambda - B) =
\lambda^4-\lambda^2+1$, which has no real zeros, so $B$ has no real
eigenvalues. Let $v = (0,0,0,1)$. We have $B^T v = Bv$, so, for every $x \in
\real^4$, we have $x \cdot Bv - v \cdot Bx = 0$.  Thus, if $h$ is any
element of~$\Lie h_B \cap \Lie n$ with $x_h = v$, then $h$ is in the center
of~$\Lie h_B \cap \Lie n$. Therefore, the center of~$\Lie h_B \cap \Lie n$
contains $\langle h, \Lie u_{\alpha+2\beta} \rangle$, so the dimension of
the center is at least~$2$. (In fact, the center is $3$-dimensional, as will
be explained below.)  Because the center of $\Lie h_{\SU} \cap \Lie n$
is~$\Lie u_{\alpha+2\beta}$, which is one-dimensional, we conclude that
$\Lie h_B$ is not isomorphic to~$\Lie h_{\SU}$. 
 \end{rem}

\begin{rem}
 Almost every $H_B$ is isomorphic to $H_{\SU}$. Namely, $H_B$ is isomorphic
to~$H_{\SU}$ if $B$ belongs to the dense, open set where $\det(B^T - B) \neq
0$. (Perturb $B$ by adding almost any skew-symmetric matrix to make the
determinant nonzero.  If the skew-symmetric matrix is small enough, the
perturbation will not have any real eigenvalues.)  To see this, note that,
for every~$B$, the torus $A \cap H_B$ has only two weights on the unipotent
radical (one on $\Lie h_B \cap (\Lie u_\beta + \Lie u_{\alpha+\beta})$, and
2~times that on $\Lie u_{\alpha+2\beta} = [\Lie h_B \cap \Lie n, \Lie h_B
\cap \Lie n]$).  So two $H_B$'s are isomorphic if and only if their
unipotent radicals are isomorphic.  We show below that each unipotent
radical is the direct product of an abelian group with a Heisenberg group,
so the unipotent radicals are isomorphic if and only if their centers have
the same dimension.  Finally, the dimension of the center of~$H_B$ is
$1+{}$the dimension of the kernel of $B^T - B$.  Therefore, it is easy to
see which $H_B$'s are isomorphic.  (Also, there are only finitely many
different $H_B$'s, up to isomorphism.)

We now show that the unipotent radical of~$H_B$ is the direct product of an
abelian group with a Heisenberg group. Let $\Lie q_0$ be the kernel of $B^T -
B$, and let $\Lie w_0$ be a subspace of $R^{2m-2}$ that is complementary
to~$\Lie q$. Define 
 $$\Lie q = \bigset{
 \begin{pmatrix}
 0 & 0 & v & 0         & 0  \\
   &0  & B(v)    & 0             &0 \\
   &   & \dots
 \end{pmatrix}
 }{v \in \Lie q_0}
 \subset \Lie h_B$$
 and
 $$\Lie w = \bigset{
 \begin{pmatrix}
 0 & 0 & w & 0         & 0  \\
   &0  & B(w)    & 0             &0 \\
   &   & \dots
 \end{pmatrix}
 }{w \in \Lie w_0}
 \subset \Lie h_B, $$
 and let $\Lie z$ be the center of~$\Lie n$. Then the
Lie algebra of the unipotent radical of $H_B$ is $(\Lie w + \Lie z) + \Lie
q$. Let us see that $\Lie w + \Lie z$ is a Heisenberg Lie algebra.  Choose
some nonzero $z_0 \in \Lie z$.  For any $v,w \in \Lie w$, there is some
scalar $\langle v\mid w \rangle$, such that
	\begin{equation}
 [v,w] = \langle v\mid w \rangle z_0 \label{Heisenberg-defn}
 \end{equation}
 (because $[\Lie w,\Lie w] \subset \Lie z$ and $\Lie z$ is
one-dimensional).  Clearly, $\langle \cdot \mid \cdot \rangle$ is skew
symmetric, because $[v,w] = -[w,v]$.  Also, for every $v \in \Lie w$, we
have $\langle v\mid \Lie w \rangle \neq 0$, because $\Lie w \cap \Lie q = 0$;
so the form is nondegenerate on~$\Lie w$.  Thus, $\langle \cdot \mid \cdot
\rangle$ is a symplectic form on~$\Lie w$, so \pref{Heisenberg-defn} is the
definition of a Heisenberg Lie algebra.
 \end{rem}

\section{Non-existence results on compact Clifford-Klein forms of
$\SO(2,n)/H$} \label{SO2n-sect}

The following lemma is obtained by combining Theorem~\ref{noncpct-dim}
with some of our results on Cartan projections.

\begin{lem} \label{noCKif<n}
  Assume that $G = \SO(2,n)$. Let $H$ be a closed, connected subgroup
of~$AN$.
 \begin{enumerate}
 \item \label{noCKif<n-linear}
 If $\mu(H) \approx \muH{\|h\|}{\cdot}$ and $\dim H < n$, then
$G/H$ does not have a compact Clifford-Klein form.
 \item \label{noCKif<n-quadratic}
 If $\mu(H) \approx \muH{\cdot}{\|h\|^2}$ and $\dim H < 2\lfloor
n/2\rfloor$, then $G/H$ does not have a compact Clifford-Klein form.
 \end{enumerate}
 \end{lem}

\begin{proof}
 By Theorem~\fullref{noncpct-dim}{-d(L)>d(H)}, it is enough to find a
subgroup~$L$ of~$AN$, such that $\dim L > \dim H$ and $L \subset CHC$ for
some compact set~$C$.

 \pref{noCKif<n-linear} We take $L$ to be $\SO(1,n) \cap AN$. We have $\dim L
= n > \dim H$. Furthermore, we know from Remark~\ref{mu(SO1n)} that $\mu(L)
\approx \muH{\|h\|}{\|h\|}$. Thus, from the assumption on the form
of~$\mu(H)$, we know that there is a compact subset~$C$ of~$G$ with $L
\subset CHC$ \see{H'inCHC}.

 \pref{noCKif<n-quadratic} Let $n' = 2 \lfloor n/2 \rfloor$. Because
$n'-2$ is even, there is a linear transformation $B \colon \real^{n'-2}
\to  \real^{n'-2}$ that has no real eigenvalues. In other words, $\dim
\langle x, Bx \rangle = 2$ for all nonzero $x \in \real^{n'-2}$. Let
 $$\Lie h_B =
 (\ker \alpha) + \{\, h \in \Lie n \mid \phi_h = 0, x_h \in \real^{n'-2},
y_h = B(x_h) \,\} ,$$
 where we identify $\real^{n'-2}$ with a codimension-one subspace
of~$\real^{n-2}$ in the case where $n$~is odd (so $n' = n-1$),  and let $H_B$
be the corresponding connected, closed subgroup of~$AN$.
 
Let $L = H_B$.  (Thus, for the appropriate choice of~$B$, we could take $L =
\SU(1,n'/2) \cap AN$.) Then  $\dim L = n' > \dim H$. We know from
Corollary~\fullref{SO2n-semi-notCDS}{-<xy>not1} that
 $\mu(L) \approx \muH{\|h\|^2}{\|h\|^2}$. Thus, from the assumption on the
form of~$\mu(H)$, we know that there is a compact subset~$C$ of~$G$ with $L
\subset CHC$ \see{H'inCHC}, as desired.
 \end{proof}

\begin{prop}[{see \ref{unip->notess}}] \label{unip->notess(cpct)}
 Assume that $G = \SO(2,n)$. If $H$ is a nontrivial, connected, unipotent
subgroup of~$G$, then $G/H$ does not have a compact Clifford-Klein form.
 \end{prop}

 \begin{notation} \label{L5-defn}
 Let $L_5 \subset \SL_5(\real)$ be the image of $\SL_2(\real)$ under an
irreducible $5$-dimensional representation of $\SL_2(\real)$. There is an
$L_5$-invariant, symmetric, bilinear form of signature $(2,3)$ on~$\real^5$,
so we may view $L_5$ as a subgroup of $\SO(2,3)$. More concretely, we may
take the Lie algebra of~$L_5$ to be the image of the homomorphism $\pi
\colon \sltwoR \to \so(2,3)$ given by
 $$
 \pi
 \begin{pmatrix}
 t & u \cr
 v & -t \\
 \end{pmatrix}
 =
 \bigset{
 \begin{pmatrix}
 4t & 2 u & 0    & 0     & 0     \cr
 2 v & 2t          & \sqrt{6} u & 0  & 0  \cr
 0  & \sqrt{6} v & 0    & -\sqrt{6} u & 0     \cr
 0  & 0    &-\sqrt{6} v & -2t    & -2 u \cr
 0          & 0    & 0    &-2 v  & -4t   \cr
 \end{pmatrix}
 }{ t, u, v \in \real}
 .
 $$
 Via the embedding $\real^5 \hookrightarrow \real^{n+2}$ given by
 $$ (x_1,x_2,x_3,x_4,x_5) \mapsto  (x_1,x_2,x_3,0,0,\ldots, 0,0,x_4,x_5) ,$$
 we may realize $\SO(2,3)$ as a subgroup of $\SO(2,n)$, so we may view $L_5$
as a subgroup of $\SO(2,n)$ (for any $n \ge 3$).
 \end{notation}

\begin{prop}[{Oh \cite[Ex.~5.6]{OhTempered}}] \label{L5-tempered}
 Assume that $G = \SO(2,3)$. Then $L_5$ is $(G,K)$-tempered
\see{tempered-def}.
 \end{prop}

\begin{cor} \label{SO/L5-notess}
 Assume that $G = \SO(2,n)$. Then $G/L_5$ does not have a compact
Clifford-Klein form.
 \end{cor}

\begin{proof}
 From Proposition~\ref{L5-tempered}, we know that $L_5$ is $\bigl(
\SO(2,3),K' \bigr)$-tempered, for any maximal compact subgroup~$K'$
of~$\SO(2,3)$. From the definition, it is clear that this implies that $L_5$
is $(G,K)$-tempered. (More generally, a tempered subgroup of any closed
subgroup of~$G$ is a tempered subgroup of~$G$.) Therefore,
Theorem~\ref{MargulisTempered} implies that $G/L_5$ does not have a compact
Clifford-Klein form.
 \end{proof}

\begin{prop}[{Kulkarni \cite[Cor.~2.10]{Kulkarni}}]
\label{SO2n/SO1odd-notess}
 If $n$~is odd, then $\SO(2,n)/SO(1,n)$
does not have a compact Clifford-Klein form.
 \end{prop}

One direction of Theorem~\ref{CK-even} in the introduction is contained in
the following theorem. The converse is obtained by combining
Theorem~\ref{compact} with Kulkarni's construction \cite[Thm.~6.1]{Kulkarni}
of a compact Clifford-Klein form of $\SO(2,n)/\SO(1,n)$ when $n$~is even. 

\begin{thm} \label{noCK}
 Assume that $G = \SO(2,n)$, with $n \ge 3$. Let $H$ be a closed, connected
subgroup of~$AN$, and such that $H$ is compatible with~$A$ \see{compatible}. 
 Assume that neither $H$ nor $G/H$ is compact, and that $G/H$ has a compact
Clifford-Klein form.
 If $n$ is even, then $H$ is one of the two types described in
Theorem~\ref{CK-even}.
 If $n$~is odd, then either
 \begin{enumerate}
 \item \label{noCK-<xy>not1}
 $\dim H = n-1$, and  $H$ is of type~\fullref{SO2n-semi-notCDS}{-<xy>not1};
or
 \item $n = 3$, $\dim H = 2$, and either
 \begin{enumerate}
 \item \label{noCK-notsemi-a,a+2b}
 $H$ is of type~\fullref{SO2n-notsemi-notCDS}{-a,a+2b}; or
 \item \label{noCK-notsemi-a+2b,a}
 $H$ is of
type~\fullref{SO2n-notsemi-notCDS}{-a+2b,a}; or
 \item \label{noCK-a,a+2b}
 $H$ is of type~\fullref{SO2n-semi-notCDS}{-root},
 and there is a positive root~$\gamma$ and a real number $p \in
(-1/3,1/3)$, such that $\{\omega,\gamma\} = \{\alpha,\alpha+2\beta\}$
and $T = \ker(p \omega + \gamma)$.
 \end{enumerate}
 \end{enumerate}
 \end{thm}

\begin{proof}
 Because of Proposition~\ref{unip->notess} and Lemma~\ref{CDS->notess}, we
know that $H \not\subset N$ and that $H$ is not a Cartan-decomposition
subgroup. Thus, $H$ must be one of the subgroups described in either
Theorem~\ref{SO2n-semi-notCDS} or Theorem~\ref{SO2n-notsemi-notCDS}. From
Proposition~\ref{tempered}, we know that $H$ cannot be of
type~\fullref{SO2n-semi-notCDS}{-1D}.

Let us first consider the cases where $\mu(H)$ is of the form
 $\mu(H) \approx \muH{\|h\|}{\cdot}$.
 \begin{itemize}
 \item If $H$ is of type~\fullref{SO2n-semi-notCDS}{-y=0}, then $\dim H
\le n-1$.
 \item If $H$ is of type~\fullref{SO2n-semi-notCDS}{-x=0}, then $\dim H
\le n-1$.
 \item If $H$ is of type~\fullref{SO2n-semi-notCDS}{-b2-2p}, then $\dim H
\le n$.
 \item If $H$ is of type~\fullref{SO2n-semi-notCDS}{-dim2}, then $\dim H
= 2$.
 \item If $H$ is of type~\fullref{SO2n-semi-notCDS}{-root} and $\omega \in
\{\beta,\alpha+\beta\}$, then $\dim H \le n-1$.
 \item If $H$ is of type~\fullref{SO2n-notsemi-notCDS}{-a,a+b},
\fullref{SO2n-notsemi-notCDS}{-a+2b,b},
\fullref{SO2n-notsemi-notCDS}{-b(a+b),w(a+2b)},
or~\fullref{SO2n-notsemi-notCDS}{-b(a+b),c(a+2b)}, then $\dim H \le n-1$.
 \end{itemize}
 Thus $\dim H < n$ except possibly when $H$ is of
type~\fullref{SO2n-semi-notCDS}{-b2-2p}. Thus, in all cases
except~\fullref{SO2n-semi-notCDS}{-b2-2p}, we conclude from
Lemma~\fullref{noCKif<n}{-linear} that $G/H$ does not have a compact
Clifford-Klein form. Now suppose that $\dim H = n$ and that $H$ is of
type~\fullref{SO2n-semi-notCDS}{-b2-2p}. From
Proposition~\ref{SO1n-whenconj} (and comparing dimensions), we see that $H$
is conjugate to $\SO(1,n) \cap (AN)$. If $n$~is even, then $H$ is listed in
Theorem~\fullref{CK-even}{-SO1n}; if $n$~is odd, then
Proposition~\ref{SO2n/SO1odd-notess} implies that $G/H$ does not have a
compact Clifford-Klein form.

We now consider the one case where $\mu(H) \approx
\muH{\|h\|^{3/2}}{\|h\|^{3/2}}$; namely, we assume that $H$ is of
type~\fullref{SO2n-semi-notCDS}{-phiy}. Then $H$ is conjugate to a
co-compact subgroup of the subgroup~$L_5$ \see{L5-defn}, so
Corollary~\ref{SO/L5-notess} implies that $G/H$ does not have a compact
Clifford-Klein form.

Finally, we consider the cases where $\mu(H)$ is of the form
 $\mu(H) \approx \muH{\cdot}{\|h\|^2}$.
 \begin{itemize}
 \item If $H$ is of type~\fullref{SO2n-semi-notCDS}{-<xy>not1}, then $\dim
H \le 2 \lfloor n/2 \rfloor$ \see{dim-SUtype}.
 \item If $H$ is of type~\fullref{SO2n-semi-notCDS}{-root} and $\omega \in
\{\alpha,\alpha+2\beta\}$, then $\dim H = 2$.
 \item If $H$ is of type~\fullref{SO2n-notsemi-notCDS}{-a,a+2b},
\fullref{SO2n-notsemi-notCDS}{-a+2b,a},
\fullref{SO2n-notsemi-notCDS}{-b(a+b),a+2b},
or~\fullref{SO2n-notsemi-notCDS}{-a+b,a}, then
$\dim H = 2$.
 \end{itemize}

Assume for the moment that $n \ge 4$. Then $2 \lfloor n/2 \rfloor > 3$,
so $\dim H < 2 \lfloor n/2 \rfloor$ except possibly when $H$ is of
type~\fullref{SO2n-semi-notCDS}{-<xy>not1}, in which case $H$ is either
listed in Theorem~\fullref{CK-even}{-SU1n} (if $n$~is even;
see~\ref{dim->HB}) or listed in Theorem~\fullref{noCK}{-<xy>not1} (if $n$~is
odd). In all other cases with $n > 3$, we conclude from
Lemma~\fullref{noCKif<n}{-quadratic} that $G/H$ does not have a compact
Clifford-Klein form.

We assume, henceforth, that $n = 3$.
 Let $H_5$ be a subgroup of~$AN$ of type~\fullref{SO2n-semi-notCDS}{-phiy}.
We know, from above, that $G/H_5$ does not have a compact Clifford-Klein
form.

 If $H$ is of
type~\fullref{SO2n-semi-notCDS}{-<xy>not1}, 
\fullref{SO2n-notsemi-notCDS}{-a,a+2b}, or
\fullref{SO2n-notsemi-notCDS}{-a+2b,a}, then $H$ is listed in
Theorem~\ref{noCK}.

 If $H$ is of type~\fullref{SO2n-notsemi-notCDS}{-b(a+b),a+2b}
or~\fullref{SO2n-notsemi-notCDS}{-a+b,a}, then there is a compact subset~$C$
of~$G$ with $H_5 \subset CHC$. Because $\dim H_5 = \dim H$, we conclude from
Lemma~\fullref{noncpct-dim}{-d(L)=d(H)} that $G/H$ does not have a compact
Clifford-Klein form.

 We may now assume that $H$ is of type~\fullref{SO2n-semi-notCDS}{-root},
and that $\omega \in \{\alpha,\alpha+2\beta\}$. If~$T$ is not of the form
described in Theorem~\fullref{noCK}{-a,a+2b}, then
Proposition~\ref{rootsemi-CDS} implies that there is a compact subset~$C$
of~$G$ with $H_3 \subset CHC$, so we conclude from
Lemma~\fullref{noncpct-dim}{-d(L)=d(H)} that $G/H$ does not have a compact
Clifford-Klein form.
 \end{proof}

We now prove a lemma used in the proof of the preceding theorem.

\begin{lem} \label{dim-SUtype}
 Assume that $G = \SO(2,n)$ and that $n$~is odd.
 If $H$ is of type~\fullref{SO2n-semi-notCDS}{-<xy>not1}, then $\dim H
\le n-1$.
 \end{lem}

\begin{proof}
 Suppose that $\dim H \ge n$. Let
 $X = \{\, x_h \mid h \in \Lie h \cap \Lie n\,\}$.
 For any $x \in X$, there is a unique $B(x) \in \real^{n-2}$ such that
there is some $h \in \Lie h \cap \Lie n$ with $x_h = x$ and $y_h = B(x)$.
(The element $B(x)$ is unique because of the assumption that $\dim\langle
x,y \rangle \neq 1$.)
 Because $n - 1 \le \dim (\Lie h\cap \Lie n) \le 1+\dim X$ (with equality
on the right if $\Lie u_{\alpha+2\beta} \subset \Lie h$) and $X \subset
\real^{n-2}$, we must have $X = \real^{n-2}$. 

Thus, $B \colon \real^{n-2} \to \real^{n-2}$ is a linear transformation, and
$\{x, Bx\}$ is linearly independent for all nonzero $x \in \real^{n-2}$. This
is impossible, because any linear transformation on a real vector space of
odd dimension has an eigenvector.
 \end{proof}

\begin{cor}[of proof] \label{dim->HB}
 Assume that $G = \SO(2,n)$, and that $n$~is even.
 If $H$ is of type~\fullref{SO2n-semi-notCDS}{-<xy>not1}, and $\dim H =
n$, then there is a linear transformation $B \colon \real^{n-2} \to
\real^{n-2}$ without any real eigenvalue, such that $H = H_B$ is the
corresponding subgroup described in Theorem~\ref{compact}.
 \end{cor}

\begin{proof}[{\bf Proof of Theorem~\ref{SU1m->complete}.}]
  Assume that $G/H$ has a compact Clifford-Klein form. Replacing $H$ by a
conjugate subgroup, we may assume that there is a closed, connected
subgroup~$H'$ of~$AN$ and a compact subgroup~$C$ of~$G$, such that $CH =
CH'$ \see{HcanbeAN}. Furthermore, we may assume that $H'$ is compatible
with~$A$ \see{conj-to-compatible}. Then $H'$ must be one of the types listed
in Theorem~\ref{noCK} (with $n = 2m+1$). Therefore, noting that $\mu \bigl(
\SU(1,m) \bigr) \approx \muH{\|h\|^2}{\|h\|^2}$ and using the calculation
of~$\mu(H')$ (see~\fullref{SO2n-semi-notCDS}{-<xy>not1},
\fullref{SO2n-notsemi-notCDS}{-a,a+2b},
\fullref{SO2n-notsemi-notCDS}{-a+2b,a}, or~\ref{rootsemi-CDS}), we see that
there is a compact subset~$C_1$ of~$G$, such that $\SU(1,m) \subset C_1 H'
C_1$. Then, because $d \bigl( \SU(1,m) \bigr) = 2m = \dim H'$, we conclude
from Theorem~\fullref{noncpct-dim}{-d(L)=d(H)} that $G/\SU(1,m)$ has a
compact Clifford-Klein form. This is a contradiction.
 \end{proof}

The following lemma is used in the above proof of
Theorem~\ref{SU1m->complete}. Although the result is known, we are unable to
locate a proof in the literature. The proof here is based on our
classification of possible compact Clifford-Klein forms
(Theorem~\ref{noCK}), but the result can also be derived from the
classification of simple Lie groups of real rank one.

\begin{lem} \label{SOodd/SU}
 Assume that $G = \SO(2,2m+1)$. There does not exist a connected,
reductive subgroup~$L$ of~$G$, such that $L$ acts properly on
$G/\SU(1,m)$, and $L \backslash G/\SU(1,m)$ is compact.
 \end{lem}

\begin{proof}
 Suppose there does exist such a subgroup~$L$. Let $L = K_L A_L N_L$ be an
Iwasawa decomposition of~$L$, and let $H = A_L N_L \subset AN$.
For any co-compact lattice~$\Gamma$ in $\SU(1,m)$, we see that the
Clifford-Klein form $\Gamma \backslash G/H$ is compact, so $H$ must be one
of the subgroups described in Theorem~\ref{noCK} (or~\ref{CK-even}). Because
$2m+1$~is odd, we know that \pref{CK-even} does not apply. 
Thus, we see, from Theorem~\ref{SO2n-semi-notCDS} (or
Proposition~\ref{rootsemi-CDS} in Case~\fullref{noCK}{-a,a+2b}), that $\mu(H)
\approx \muH{\cdot}{\|h\|^2}$. Because $\mu \bigl( \SU(1,m) \bigr) \approx
\muH{\|h\|^2}{\|h\|^2}$, we conclude (e.g., from~\ref{H'inCHC}) that $H$
does not act properly on $G/\SU(1,m)$. This contradicts the fact that $L$
acts properly on $G/\SU(1,m)$.
 \end{proof}

\section{Finite-volume Clifford-Klein forms} \label{finvol-sect}

\begin{defn}[{cf.~\cite[Def.~2.2]{Kobayashi-properaction}}]
 Let $H$ be a closed, connected subgroup of~$G$, such that $G/H$ has a
$G$-invariant regular Borel measure. (Because $G$ is unimodular, this means
that $H$ is unimodular \cite[Lem.~1.4, p.~18]{Raghunathan}.)
We say that $G/H$ has a \emph{finite-volume Clifford-Klein form} if there is
a discrete subgroup~$\Gamma$ of~$G$, such that 
 \begin{itemize}
 \item $\Gamma$ acts properly on $G/H$; and 
 \item there is a Borel subset~$\mathcal{F}$ of $G/H$, such that
$\mathcal{F}$ has finite measure, and $\Gamma \mathcal{F} = G/H$.
 \end{itemize}
 \end{defn}

Unfortunately, the study of finite-volume Clifford-Klein forms does not
usually reduce to the case where $H \subset AN$, because the subgroup~$H'$ of
Proposition~\ref{HcanbeAN} is usually not unimodular.

\begin{thm} \label{SO2nCK-finvol}
 Assume that $G = \SO(2,n)$. Let $H$ be a closed, connected subgroup
of~$G$. If $G/H$ has a finite-volume Clifford-Klein form, then either
 \begin{enumerate}
 \item $H$ has a co-compact, normal subgroup that is conjugate under
$\On(2,n)$ to the identity component of either $\SO(1,n)$, $\SU(1,\lfloor
n/2 \rfloor)$, or~$L_5$ \see{L5-defn}; or
 \item \label{SO2nCK-finvol-1D}
 $d(H) \le 1$ \see{d(H)-defn}; or
 \item $H = G$.
 \end{enumerate}
 \end{thm}

It is natural to conjecture that $\SO(2,2m+1)/\SU(1,m)$ and $\SO(2,n)/L_5$ do
not have finite-volume Clifford-Klein forms, and that $G/H$ does not have a
finite-volume Clifford-Klein form when $d(H) = 1$, either.

To prepare for the proof of Theorem~\ref{SO2nCK-finvol}, we present some
preliminary results.

Unfortunately, we do not have an analogue of
Theorem~\fullref{noncpct-dim}{-d(L)>d(H)} for finite-volume Clifford-Klein
forms. The following lemma is a weak substitute.

\begin{lem} \label{HinH'->notess(finvol)}
 Let $H$ be a closed, connected, unimodular subgroup of~$G$. Suppose that
there is a closed subgroup $L$ of~$G$ containing $H$, such that $L \subset
CHC$, for some compact subset~$C$. 
 If $L/H$ does not have finite {\upshape(}$L$-invariant{\upshape)} volume,
then $G/H$ does not have a finite-volume Clifford-Klein form.
 \end{lem}

\begin{proof}
 Suppose that $\Gamma \backslash G/H$ is a finite-volume Clifford-Klein form
of $G/H$. Because $L \subset CHC$, we know that $\Gamma$ is proper on $G/L$.
(In particular, $\Gamma \cap L$ must be finite.) We have a quotient map
$\Gamma \backslash G/H \to \Gamma \backslash G/L$, and, because $\Gamma
\backslash G/H$ has finite volume, (almost) every fiber has finite measure
(cf.~\cite[\S3, p.~26--33]{Rohlin}). Thus (perhaps after
replacing $\Gamma$ by a conjugate subgroup), $(\Gamma \cap L) \backslash
L/H$ has finite volume. Because $\Gamma \cap L$ is finite, this implies that
$L/H$ has finite volume, which is a contradiction.
 \end{proof}

We now state two special cases of the Borel Density Theorem.

\begin{lem}[{cf.~\cite[Thm.~5.5, p.~79]{Raghunathan}}] \label{BDT}
 Let $H$ be a closed subgroup of~$G$. If $H$ is connected, and  $G/H$ has
finite {\upshape(}$G$-invariant{\upshape)} volume, then $H = G$.
 \end{lem}

\begin{lem}[{Mostow \cite[Prop.~11.2]{MostowFSS}, \cite[Thm.~3.1,
p.~43]{Raghunathan}}] \label{BDT-solv} 
 Suppose that $R$ is a simply connected, solvable Lie group. If $H$ is a
closed, connected subgroup of~$R$, such that $R/H$ either is compact or has
finite {\upshape(}$R$-invariant{\upshape)} volume, then $H = R$.
 \end{lem}

The following lemma is analogous to the fact
\cite[Rmk.~1.9, p.~21]{Raghunathan} that a locally compact group that
admits a lattice must be unimodular. 

\begin{lem}[{cf.~\cite[Prop.~2.1]{Zimmer}}] \label{normalizer}
 Let $H$ be a unimodular subgroup of~$G$, and
assume that there is an element~$a$ of the normalizer of~$H$ such that
the action of~$a$ by conjugation on~$H$ does not preserve the Haar
measure on~$H$. Then $G/H$ has neither compact nor
finite-volume Clifford-Klein forms.
 \end{lem}

\begin{prop} \label{unip->notess}
 Assume that $G = \SO(2,n)$. If $H$ is a nontrivial, connected, unipotent
subgroup of~$G$, then $G/H$ has neither a compact nor a finite-volume
Clifford-Klein form.
 \end{prop}

\begin{proof}
 Replacing $H$ by a conjugate, we may assume that $H$ is contained in~$N$.
Furthermore, because Lemma~\ref{BDT} implies that $G/H$ does not have finite
volume, we see from Lemma~\ref{HinH'->notess(finvol)} that we may assume
that $H$ is not a Cartan-decomposition subgroup. The proof now breaks up
into cases, determined by Proposition~\ref{HinN}. Because $H$ is unimodular,
any compact Clifford-Klein form would also have finite volume, so we need
only consider the more general finite-volume case.

\setcounter{case}{0}

\begin{case}
 Assume that $\dim H \le 1$.
 \end{case}
 Because $H$ is nontrivial and connected, we must have $\dim H = 1$.
Then $H$ is contained in a subgroup~$S$ of~$G$ that is locally
isomorphic to~$\SL(2,\real)$, and we have $H = N_S$, where  $S = K_S A_S
N_S$ is an Iwasawa decomposition of~$S$  \cite[Thm.~17(1),
p.~100]{Jacobson}. The subgroup $A_S$ normalizes~$H$, but does not
preserve the Haar measure on~$H$, so we conclude from
Lemma~\ref{normalizer} that $G/H$ does not have a finite-volume
Clifford-Klein form.

Henceforth, we assume that $\dim X \ge 2$.

\begin{case}
 Assume, for every nonzero element~$h$ of~$\Lie h$, that we have $\phi_h = 0$
and $\dim \langle x_h,y_h \rangle \neq 1$.
 \end{case}
 From Lemma~\fullref{HinN}{-<xy>not1}, we have $\mu(H) \approx
\muH{\|h\|^2}{\|h\|^2} \approx \mu(U_{\alpha+2\beta} H)$. Thus, if $G/H$ has
a finite-volume Clifford-Klein form, then Lemma~\ref{HinH'->notess(finvol)}
implies that $U_{\alpha + 2\beta}H/H$ must have finite volume. Therefore,
Lemma~\ref{BDT-solv} implies that $U_{\alpha + 2\beta}H/H$ must be
trivial, so $U_{\alpha + 2\beta} \subset H$. Then $\ker\alpha$
normalizes~$H$, but the action of $\ker\alpha$ by conjugation on~$H$ does
not preserve the Haar measure on~$H$, so Lemma~\ref{normalizer} implies that
$G/H$ does not have a finite-volume Clifford-Klein form.

\begin{case}
 Assume, for every nonzero element~$h$ of~$\Lie h$, that we have $\phi_h = 0$
and $\dim \langle x_h,y_h \rangle = 1$.
 \end{case}
 Because $\dim \Lie h \ge 2$, we know that $\Lie h \not\subset \Lie
u_{\alpha+2\beta}$, so it cannot be the case that both of $x_h$ and~$y_h$
are~$0$, for every $h \in \Lie h$. Therefore, by perhaps
replacing~$H$ with its conjugate under the Weyl reflection corresponding to
the root~$\alpha$, we may assume that $x_{h_0} \neq 0$ for some  $h_0 \in
\Lie h$. Then, because $\dim\langle x_h,y_h \rangle = 1$ for every $h \in
\Lie h$, it follows that there is a real number~$p$, such that for all $h
\in \Lie h$, we have $y_h = px_h$. Therefore, by replacing~$H$ with a
conjugate under~$U_{-\alpha}$, we may assume that $y_h = 0$ for every $h \in
\Lie h$. So $\Lie h \subset \Lie u_{\alpha+\beta} + \Lie u_{\alpha +
2\beta}$, so $\ker\beta$ normalizes~$H$. But the action of $\ker\beta$ by
conjugation on~$H$ does not preserve the Haar measure on~$H$, so
Lemma~\ref{normalizer} implies that $G/H$ does not have a finite-volume
Clifford-Klein form.

\begin{case}
 Assume that there exists a subspace~$X_0$ of~$\real^{n-2}$, $b \in X_0$,
$c \in X_0^\perp$, and $p \in \real$ with $\|b\|^2 - \|c\|^2 - 2p < 0$,
such that for every element of~$\Lie h$, we have $y= 0$, $x \in \phi c +
X_0$, and $\eta = p\phi + b \cdot x$.
 \end{case}
 Replacing $H$ by a conjugate, we may assume that $H$ is contained in
$G' = \SO(1,n)$ \see{SO1n-whenconj}. Then $H$ is a Cartan-decomposition
subgroup of~$G'$ \see{Rrank1-CDS}, but Lemma~\ref{BDT} implies that $G'/H$
does not have finite volume, so we conclude from
Lemma~\ref{HinH'->notess(finvol)} that $G/H$ does not have a finite-volume
Clifford-Klein form.
 \end{proof}

\begin{proof}[{\bf Proof of Theorem~\ref{SO2nCK-finvol}}]
 By replacing $H$ with a conjugate, we may assume that there is a compact
subgroup~$C$ of~$G$ and a closed, connected subgroup~$H'$ of~$AN$, such that
$CH = CH'$ \see{HcanbeAN}. (Note that if $H'$ is unimodular, then it is easy
to see that $G/H'$, like $G/H$, has a finite-volume Clifford-Klein form.) We
may assume that $H$ is not a Cartan-decomposition subgroup
(see~\ref{HinH'->notess(finvol)} and~\ref{BDT}), so $H'$ is not a
Cartan-decomposition subgroup.  From Proposition~\ref{unip->notess}, we see
that $H' \not\subset N$. Then, assuming, as we may, that $H'$ is compatible
with~$A$ \see{conj-to-compatible}, the subgroup $H'$ must be one of the
subgroups described in either Corollary~\ref{SO2n-semi-notCDS} or
Corollary~\ref{SO2n-notsemi-notCDS}. Furthermore, we may assume that $\dim
H' \ge 2$, for otherwise, conclusion~\pref{SO2nCK-finvol-1D} holds.

The only unimodular subgroups listed in either
Corollary~\ref{SO2n-semi-notCDS} or Corollary~\ref{SO2n-notsemi-notCDS} are
the subgroups of type~\fullref{SO2n-semi-notCDS}{-1D}, which are
one-dimensional. Therefore, $H'$ is not unimodular. Because $CH' = CH$ is
unimodular, this implies that $C$ does not normalize~$H'$, so we see
from~\ref{HcanbeAN} that $H/\Rad H$ is not compact. In other words, we may
assume that there is a connected, noncompact, semisimple Lie group~$L$ with
no compact factors, and a closed, connected subgroup~$R$ of~$AN$, such that
$H = L \semiprod R$. Because $H$ is not a Cartan-decomposition subgroup, we
know that $A \not\subset L$, so $\Rrank L = 1$. Thus, $L \cap A$ is
one-dimensional, so $H' \cap A$ is nontrivial. Thus, $H'$ cannot be any of
the subgroups listed in Corollary~\ref{SO2n-notsemi-notCDS}. In other words,
we have $H' = (H' \cap A) \semiprod (H' \cap N)$. Also, because $A
\not\subset H$, we know that that $R \subset N$ \see{HN=AN}. 

\setcounter{case}{0}

\begin{case}
 Assume that $R$ is nontrivial.
 \end{case}
 Because $R$ is a unipotent, normal subgroup of~$H$, a fundamental
theorem of Borel and Tits \cite[Prop.~3.1]{BorelTits-unip} implies
that there is a parabolic subgroup~$P$ of~$G$, such that $H$ is contained
in~$P$, and $R$ is contained in the unipotent radical of~$P$. Now, because
$\Rrank G = 2$, there are, up to conjugacy, only two parabolic subgroups
of~$G$ whose maximal connected semisimple subgroups have real rank
one.

\begin{subcase}
 Assume that $L = \langle U_\alpha, U_{-\alpha} \rangle$ and that $R
\subset U_\beta U_{\alpha+\beta} U_{\alpha+2\beta}$.
 \end{subcase}
 Every nontrivial $L$-invariant subalgebra of $\Lie u_\beta +
\Lie u_{\alpha+\beta} + \Lie u_{\alpha+2\beta}$ contains~$\Lie
u_{\alpha+2\beta}$, so $\Lie r$ must contain~$\Lie
u_{\alpha+2\beta}$. Then $\ker \alpha$ normalizes~$L$ and~$R$ but does not
preserve the Haar measure on~$LR$, so Lemma~\ref{normalizer} implies that
$G/H$ does not have a finite-volume Clifford-Klein form.

\begin{subcase}
 Assume that $L = \langle U_\beta, U_{-\beta} \rangle$ and that $R
\subset U_\alpha U_{\alpha+\beta} U_{\alpha+2\beta}$.
 \end{subcase}
 Then $\ker \beta$ normalizes~$L$ and~$R$, but does not preserve the
Haar measure on~$LR$, so Lemma~\ref{normalizer} implies that $G/H$
does not have a finite-volume Clifford-Klein form.

\begin{case}
 Assume that $R$ is trivial.
 \end{case}
 This means that $H$ is reductive. Then, because $\Rrank H = 1$, we may
assume that $H$ is almost simple, by removing the compact factors of~$H$.
Thus, we see from Lemma~\ref{SO2n-maxsimple}
(and~\ref{HinH'->notess(finvol)}) that the identity component of~$H$ is
conjugate under $\On(2,n)$ to the identity component of either $\SO(1,n)$,
$\SU(1,\lfloor n/2 \rfloor)$, or~$L_5$, as desired.
 \end{proof}

We now prove a classification theorem we used in the above proof. Although
this result is known, the authors are not aware of any convenient reference.

\begin{lem} \label{SO2n-maxsimple}
 Assume that $G = \On(2,n)$. If $L$ is any connected, almost-simple
subgroup of~$G$, such that $\Rrank L = 1$, then $L$ is conjugate to a
subgroup of either $\SO(1,n)$, $\SU(1,\lfloor n/2 \rfloor)$, or~$L_5$.
 \end{lem}

\begin{proof}[Sketch of proof]
 Let $L = K_L A_L N_L$ be an Iwasawa decomposition of~$L$, and let $H = A_L
N_L$. We may assume that $A_L \subset A$ and that $N_L \subset N$. Because
$\mu(L) \approx \mu(A_L) \not\approx A^+$, we know that $H$ is not a
Cartan-decomposition subgroup of~$G$, so $H$ must be one of the subgroups
described in Theorem~\ref{SO2n-semi-notCDS}. Because $H$ is an epimorphic
subgroup of~$L$ \cite[\S2]{BienBorel}, it suffices to show that $H$ is
conjugate to a subgroup of either $\SO(1,n)$, $\SU(1,\lfloor n/2 \rfloor)$,
or~$L_5$.

Because $N_L$ is nontrivial, we know that $H$ is not of
type~\fullref{SO2n-semi-notCDS}{-1D}.

If $H$ is of type~\fullref{SO2n-semi-notCDS}{-y=0},
\fullref{SO2n-semi-notCDS}{-x=0}, \fullref{SO2n-semi-notCDS}{-b2-2p},
or~\fullref{SO2n-semi-notCDS}{-dim2}, then $H$ is conjugate to a subgroup of
$\SO(1,n)$ \cf{SO1n-whenconj}.

If $H$ is of type~\fullref{SO2n-semi-notCDS}{-phiy}, then $H$ is conjugate
to a subgroup of~$L_5$.

If $H$ is of type~\fullref{SO2n-semi-notCDS}{-root}, then, because $\mu(H)
\approx \mu(A_L)$, we see from Proposition~\ref{rootsemi-CDS} that $p = 0$.
Replacing $H$ by a conjugate under the Weyl group, we may assume that
$\omega \in \{ \alpha+\beta, \alpha+2\beta \}$. Therefore, $H$ is contained
in either $\SO(1,n)$ or $\SU(1,\lfloor n/2 \rfloor)$.

We may now assume that $H$ is of type~\fullref{SO2n-semi-notCDS}{-<xy>not1}.

\setcounter{case}{0}

\begin{case} 
 Assume that $\Lie u_{\alpha+2\beta} \subset \Lie n_L$.
 \end{case}
  There is a subspace~$V$ of~$\real^{n-2}$ and a linear transformation $B
\colon V \to \real^{n-2}$, such that
 $$ \Lie n_L = \{\, h \in \Lie n \mid x_h \in V, \, y_h = B x_h, \, \phi_h =
0 \,\} .$$ 
 Because $[\Lie n_L,\Lie n_L] = \Lie u_{\alpha+2\beta}$  is one-dimensional,
the classification of simple Lie groups of real rank one
($\SO(1,k)$, $\SU(1,k)$, $\Sp(1,k)$, $F_4^{-20}$) implies that $L$
is locally isomorphic to $\SU(1,k)$, where $k-1 = (\dim V)/2$.

 For any $h_0 \in \Lie n_L \setminus \Lie u_{\alpha+2\beta}$, we have 
 $\langle h_0, \Lie u_{\alpha+2\beta}, \Lie u_{-(\alpha+2\beta)} \rangle \iso
\Lie{su}(1,2)$. In particular, if $h_0$ is the element of $\SO(2,4)$ with
 $x_{h_0} = (1,0)$ and $y_{h_0} = (0,1)$, then
 $$ \langle h_0, \Lie u_{\alpha+2\beta}, \Lie u_{-(\alpha+2\beta)} \rangle
\cap \Lie n = 
  \bigset{
  \begin{pmatrix}
 0 & 0 & s & t  & \eta &0\\
   & 0 & -t& s  &   0 & -\eta \\
   &   & \cdots &  \\
 \end{pmatrix}
 }{ s,t, \eta \in \real}
 .$$
  From this, we conclude that $B(V) = V$. 

The desired conclusion is easy if $k = 1$, so let us assume that $k \ge 2$.
From the structure of $\SU(1,k)$, we know that $C_L(A_L)$ contains a
subgroup~$M$ that is isomorphic to $\SU(k-1)$ and acts transitively on the
unit sphere in $\Lie n_L/\Lie u_{\alpha+2\beta}$. Then the action of~$M$
on~$V$ is transitive on the unit sphere in~$V$ and, because $M$ normalizes
$\Lie n_L$, we see that $B$ is $M$-equivariant. Therefore $B$ is a scalar
multiple of an orthogonal transformation, so $H$ is conjugate to a subgroup
of $\SU(1, \lfloor n/2 \rfloor)$. 

\begin{case} 
 Assume that $\Lie u_{\alpha+2\beta} \not\subset \Lie n_L$.
 \end{case}
 There is a subspace~$V$ of~$\real^{n-2}$ and a linear transformation $B
\colon V \to \real^{n-2}$, such that
 $$ \Lie n_L = \{\, h \in \Lie n \mid x_h \in V, \, y_h = B x_h, \, \phi_h =
\eta_h = 0 \,\} .$$ 
 Note that $N_L$ must be abelian (because $\Lie n_L \cap \Lie
u_{\alpha+2\beta} = 0$), so, from the classification of simple Lie groups of
real rank one, we see that $L$ is locally isomorphic to $\SO(1,k)$, where
$k-1 = \dim V$. 

Assume, for the moment, that $k \ge 4$, so $\SO(k-1)$ is almost simple. 
From the structure of $\SO(1,k)$, we see that $C_L(A_L)$ contains a
subgroup~$M$ that is isomorphic to $\SO(k-1)$ and acts transitively on the
unit spheres in~$V$ and~$B(V)$, such that $B$ is $M$-equivariant. (In
particular, the ratio $\|Bv\|/\|v\|$ is constant on $V \setminus \{0\}$.)
Therefore, letting $\pi \colon \real^{n-2} \to V$ be the orthogonal
projection, we see that the composition $\pi \circ B \colon V \to V$ must be
a real scalar. Replacing $L$ by a conjugate under~$U_{-\alpha}$, we may
assume that this scalar is~$0$; thus, $B(V)$ is orthogonal to~$V$, so $H$ is
conjugate to a subgroup of $\SU(1, \lfloor n/2 \rfloor)$, as desired.

Now assume that $k<4$. The desired conclusion is easy if
$k=2$, so we may assume that $k = 3$. Let $V_1$ be an irreducible summand of
the $L$-representation on $\real^{n-2}$, and let $\langle\mid\rangle$ be the
$\SO(2,n)$-invariant bilinear form on~$\real^{n+2}$. 

\begin{subcase}
 Assume that the restriction of $\langle\mid\rangle$ to~$V_1$ is~$0$.
 \end{subcase}
 There must be $L$-invariant subspaces $V_0$ and~$V_{-1}$, such that
$\real^{n+2} = V_{-1} \oplus V_0 \oplus V_1$, and $\langle V_i \mid V_j
\rangle = 0$ unless $i = -j$. Thus, the restriction of $\langle\mid\rangle$
to $V_{-1} \oplus V_1$ is split. Because $\langle\mid\rangle$ has signature
$(2,n)$, we conclude that $\dim V_{\pm 1} \le 2$. Because the smallest
nontrivial representation of~$L$ is 3-dimensional, this implies that $L$
acts trivially on $V_{-1} \oplus V_1$. Therefore, $L$ fixes a vector of
norm~$-1$, so $L$ is contained in a conjugate of $\SO(1,n)$.

\begin{subcase} \label{simplesubs-subcase1m}
 Assume that the restriction of $\langle\mid\rangle$ to~$V_1$ has signature
$(1,m)$, for some~$m$.
 \end{subcase}
 Let $V_1 = W_{-r} \oplus \cdots \oplus W_{r-1} \oplus W_r$ be the
decomposition of~$V_1$ into weight spaces (with respect to~$A_L$). We must
have $\langle W_i \mid W_j \rangle = 0$ unless $i = -j$. From the assumption
of this subcase, we conclude that $r = 1$ and $\dim W_1 = 1$. It follows
that $V_1$ is the standard representation of~$L$ on~$\real^4$ (whose
$L$-invariant bilinear form is unique up to a real scalar). 

There must be another irreducible summand~$V_2$, such that the restriction of
$\langle\mid\rangle$ to~$V_2$ has signature $(1,*)$, and $\real^{n+2} = V_1
\oplus V_2$. From the argument of the preceding paragraph, we see that $\dim
V_2 = 4$, and that the representations of~$L$ on $V_1$ and~$V_2$ are
isomorphic (so $L$ acts diagonally on $\real^4 \oplus \real^4 \iso
\real^8$). Then $L$ is conjugate to a subgroup of $\SU(1,3)$.

\begin{subcase}
 Assume that $V_1 = \real^{n+2}$.
 \end{subcase}
  Let $V_1 = W_{-r} \oplus \cdots \oplus W_{r-1} \oplus W_r$ be the
decomposition of~$V_1$ into weight spaces (with respect to~$A_L$). 

If $\dim W_r = 2$, then $V_1$ is an irreducible $\complex$-representation of
$\so(1,3) \iso \sltwoC$, so $\dim W_i = 2$ for
all~$i$. From the signature of $\langle\mid\rangle$, we must have $r = 1$,
so $V_1$ is the 3-dimensional irreducible $\complex$-representation of
$\sltwoC$, that is, the adjoint representation.
However, the invariant bilinear form $\operatorname{Re}\bigl(
\operatorname{Tr} X^2 \bigr)$ for this representation has signature
$(3,3)$, not $(2,n)$. This is a contradiction.

We now know that $\dim W_r = 1$. Also, we must have $r = 2$ (because the
argument of Subcase~\ref{simplesubs-subcase1m} applies if $r=1$). The
signature of $\langle\mid\rangle$ implies that $\dim W_1 = \dim W_2 = 1$.
There is no such irreducible representation of $\so(1,3)$, so this is a
contradiction.
 \end{proof}

\section{Non-existence of compact Clifford-Klein forms of $\SL(3, \real)/H$}
\label{SL3-sect}

Y.~Benoist \cite[Cor.~1]{Benoist} proved that $\SL(3,\real)/\SL(2,\real)$
does not have a compact Clifford-Klein form. By combining Benoist's method
with the fact that if $G/H$ has a compact Clifford-Klein form, then $H$
cannot be one-dimensional \see{tempered}, we show, more generally, that
no interesting homogeneous space of $\SL(3,\real)$ has a compact
Clifford-Klein form.

\begin {thm}[{Benoist \cite[Thms.~3.3 and 4.1]{Benoist} (see
also Prop.~\ref{Benoist-NinH})}] \label{Ben-B}
 Let $H$ be a closed, connected subgroup of~$G$, such that, for some compact
set~$C$ in~$A$, we have $B^+ \subset \mu (H) C$, where $B^+$ is defined
in Notation~\ref{oppinv}. Then $G/H$ has neither compact
nor finite-volume Clifford-Klein forms, unless $G/H$ is compact.
 \end{thm}

\begin{notation} \label{oppinv}
 Let $i$ be the opposition involution in~$A^+$, that is, for $a\in A^+$,
$i(a)$ is the unique element of~$A^+$ that is conjugate to $a^{-1}$, and
set $B^+=\{\,a\in A^+ \mid i(a)=a \,\}$.
 \end{notation}

\begin{prop} \label{SL3-B+}
 Assume that $G = \SL(3,\real)$.
 Let $H$ be a closed, connected subgroup of~$G$ with $d(H) \ge 2$
\see{d(H)-defn}. Then $B^+ \subset \mu(H)$.
 \end{prop}

\begin{proof}
 Since $H$ is non-compact, there is a curve $h_t$ in $H$ such that $h_0=e$
and $h_t \to \infty$ as $t\to \infty$.  Since $d(H) \ge 2$ (so
$H/K_H$ is homeomorphic to~$\real^k$, for some $k \ge 2$), it is easy to find
a continuous and proper map $\Phi\colon [0,1] \times \real ^+ \to H$ such
that $\Phi (0,t)= h_t$ and $\Phi(1,t) = h_t^{-1}$, for all $t \in \real^+$. 

 If we identify the Lie algebra~$\Lie a$ of~$A$ with the connected component
of~$A$ containing~$e$, then~$A^+$ is a convex cone in~$\Lie a$ and the
opposition involution~$i$ is the reflection in~$A^+$ across the ray~$B^+$.
Thus, for any $a\in A^+$, the points $a$ and~$i(a)$ are on opposite sides
of~$B^+$, so any continuous curve in~$A^+$ from~$a$ to~$i(a)$ must
intersect~$B^+$. In particular, any curve from~$\mu(h_t)$ to
$\mu(h_t^{-1})$ must intersect~$B^+$. Thus, we see, from an
elementary continuity argument, that $\mu \bigl[ \Phi\bigl( [0,1] \times
\real ^+ \bigr) \bigr]$ contains~$B^+$. Therefore, $B^+$ is contained in
$\mu(H)$.
 \end{proof}

\begin{proof}[{\bf Proof of Proposition~\ref{no-cpt}}]
 Suppose that $H$ is non-compact. By Proposition~\ref{tempered}, we may
assume that $d(H) \ge 2$. Then, from Proposition~\ref{SL3-B+}, we know that
$\mu(H) \supset B^+$, so Theorem~\ref{Ben-B} implies that $G/H$ is compact. 
 \end{proof}

The following is proved similarly.

\begin{cor} \label{no-finvol}
 Let $H$ be a closed, connected subgroup of $G = \SL(3,\real)$. If $G/H$
has a finite-volume Clifford-Klein form, then either $d(H) \le 1$, or $H =
G$.
 \end{cor}

\begin{appendix}

\section{A short proof of a theorem of Benoist}
\label{Benoist-sect}

Most of Benoist's paper~\cite{Benoist} is stunningly elegant; the only
exception is Section~4, which presents a somewhat lengthy argument to
eliminate a troublesome case. We provide an alternative treatment of this
one case. Our proof does not match the elegance of the rest of Benoist's
paper, but it does get the unpleasantness over fairly quickly.

Our version of the result is not a complete replacement for Benoist's,
because we require the action of $\Gamma$ on $G/H$ to be proper, while
Benoist makes no such assumption. There are many situations where one is
interested in improper actions (for example, a quotient of a proper action
is usually not proper), but, in applications to Clifford-Klein forms, the
action is indeed proper, so our weaker version of the proposition does apply
in that situation. Our proof has the virtue that, for the real field, it
does not require the subgroup~$H$ to be Zariski closed. It also applies to
finite-volume Clifford-Klein forms, not just compact ones.

\begin{prop}[{Benoist \cite[Thm.~4.1]{Benoist}}] \label{Benoist-NinH}
 Let $\G$ be a Zariski connected, semi\-simple, algebraic group over a local
field~$k$ of characteristic~$0$, and let $G = \G_k$. Suppose that
$\Gamma$ and~$H$ are closed subgroups of~$G$, and assume that
 \begin{enumerate}
 \item $\Gamma$ is nilpotent;
 \item $\Gamma$ acts properly on~$G/H$; 
 \item either
 \begin{enumerate}
 \item $\Gamma\backslash G/H$ is compact, or
 \item $\Gamma$ is discrete and $\Gamma\backslash G/H$ has finite volume; and
 \end{enumerate}
 \item  either
 \begin{enumerate}
 \item $H$ is Zariski closed, or
 \item $k = \real$ or~$\complex$, and $H$ is almost connected.
 \end{enumerate}
 \end{enumerate}
 Then $H$ contains a conjugate of~$N$.
 \end{prop}

\begin{proof}
 To clarify the exposition, let us assume that $\Gamma$ is abelian, that
$H = (H \cap A) \semiprod (H \cap N)$, and that $\Gamma\backslash G/H$
is compact. Remark~\ref{Benoist-nilp} describes appropriate
modifications of the proof to eliminate these assumptions.

 Let $X$ be the Zariski closure of~$\Gamma$. By replacing $\Gamma$ with a
finite-index subgroup, we may assume that $X$ is Zariski connected. Note
that $X$, like~$\Gamma$, is abelian. Then, by replacing $\Gamma$ with a
conjugate subgroup, we may assume that $X = (X \cap A) (X \cap N) E$, where
$E$ is an anisotropic torus, hence compact. Because we may replace $\Gamma$
with $(\Gamma E) \cap (AN)$, there is no harm in assuming that $\Gamma
\subset AN$. 

\begin{notation}[and remarks]
 \ 
 \begin{itemize}

 \item We define a right-invariant metric on~$G$ by $d(g,h) = \log \| g
h^{-1} \|$.

    \item Let $\pi  \colon AN \to A$ and $\nu \colon AN
\to N$ be the canonical projections, so $g = \pi(g) \nu(g)$, for all
$g \in AN$. Note that $\pi$ is a homomorphism, but $\nu$ is not. 
Also note that $\pi(H) = H \cap A$ and $\nu(H) = H \cap N$, because
$H = (H \cap A) \semiprod (H \cap N)$. Also note that $\pi(\Gamma)
\subset X \cap A$ and $\nu(\Gamma) \subset X \cap N$, so
$\pi(\Gamma)$ and $\nu(\Gamma)$ centralize each other. Thus, the
restriction of~$\nu$ to~$\Gamma$ \emph{is} a homomorphism.

    \item Because $\pi(\Gamma) \subset A$, and the centralizer of any
Zariski-connected subgroup of~$A$ is a Zariski-connected, reductive
$k$-subgroup of~$G$ \cite[Thm.~22.3, p.~140]{Humphreys}, we
may write $C_G \bigl( \pi(\Gamma) \bigr ) = LZ$, where $L$ is a
Zariski-connected, reductive $k$-subgroup of~$G$ with compact center, and $Z$
is a Zariski closed (abelian) subgroup of~$A$ that centralizes~$L$. We have
$\Gamma \subset LZ$ and $\pi(\Gamma) \subset Z$. Note that we must have
$\nu(\Gamma) \subset L$, because $Z$, being a subgroup of~$A$, has no
nontrivial unipotent elements. 

    \item Let $\mul \colon L \to A \cap L$ be the Cartan projection for the
reductive group~$L$. We define $\muLZ \colon LZ \to A$ by
$\muLZ(lz) = \mul(l)z$, for $l \in L$ and $z \in Z$. (Unfortunately, $\muLZ$
is not well defined if $L \cap Z \neq e$. On the other hand, at key points of
the proof, we only calculate $\muLZ$ up to bounded error, so, because $L \cap
Z$ is finite, this is not an important issue.) In particular, for $\gamma \in
\Gamma$, we have $\muLZ(\gamma) = \mul\bigl( \nu(\gamma) \bigr)
\pi(\gamma)$. Note that, because $\muLZ(g) \in (K \cap L) g (K \cap L)$, we
have $\mu(g) = \mu \bigl( \muLZ(g) \bigr)$, for every $g \in LZ$. Therefore,
for any $g \in LZ$ and $a \in A$, we have $d \bigl( \mu(g), \mu(a)) \le d
\bigl( \muLZ(g), a \bigr)$. Thus, if we find a sequence $\gamma_n \to
\infty$ in~$\Gamma$, with $d\bigl(\muLZ(\gamma_n), H \cap A \bigr) = O(1)$,
then we have obtained a contradiction to the fact that $\Gamma$ acts
properly on $G/H$.
 \end{itemize}
 \end{notation}

Because $\Gamma\backslash AN/H$ is a closed subset of the compact
space $\Gamma\backslash G/H$, we know that it is compact, so we see,
by modding out~$N$, that $A/\closure{\pi\bigl(\Gamma H \bigr)}$ is
compact, so there is a free abelian subgroup~$\GamPrime$ of~$\Gamma$, such
that 
 \begin{enumerate}
 \item \label{sh-latt}
 $\pi(\GamPrime)$ is a co-compact, discrete subgroup of $A/\pi(H)$; and
 \item \label{sh-inj}
 $\GamPrime \cap (NH) = e$.
 \end{enumerate}
 We may assume that $G/H$ is not compact (else $AN/H$ is compact, so
Lemma~\ref{BDT-solv} implies that $H \supset N$, as desired). Then $H$
cannot be a Cartan-decomposition subgroup of~$G$, so $\dim \bigl(\pi(H)
\bigr) < \dim A$ \see{HN=AN}. Therefore, $A/\pi(H)$ is not compact, so
$\GamPrime$ is infinite.

Suppose, for the moment, that $\GamPrime$ is co-compact in~$\Gamma$. (In
this case, there is no harm in assuming that $\Gamma = \GamPrime$.) Then,
from~\pref{sh-latt}, we see that $\pi(\Gamma H)$ is closed, so the
inverse image in $\Gamma \backslash AN /H$ must be compact; that is,
$\Gamma \backslash \Gamma NH/H$ is compact. Therefore,
from~\pref{sh-inj}, we see that $NH/H$ is compact. Then $N/(H \cap N)$,
being homeomorphic to $NH/H$, is compact. Therefore, $N \subset H$
\see{BDT-solv} as desired.

We may now assume that $\Gamma/\GamPrime$ is not compact. Then, for any
$R > 0$, there is some $\gamma_0 \in \Gamma$, such that
$d(\gamma_0, \GamPrime) > R$.  (Let $\gamma_0$ be any element of~$\Gamma$
that is not $B_R(e) \GamPrime$, where $B_R(e)$ is the closed ball of
radius~$R$ around~$e$.) Because $\pi(\GamPrime)$ is co-compact in $A/\pi(H)$,
there is some $\gamprime \in \GamPrime$, such that $d\bigl(\muLZ(\gamma_0)
\pi(\gamprime) , \pi(H) \bigr) < C$, where $C > 0$ is an appropriate
constant that is independent of~$R$, $\gamma_0$, and~$\gamprime$.

Now, to clarify the argument, assume, for the moment, that $\GamPrime
\subset Z$. In this case, we see from the definition that
$\muLZ(\gamma_0 \gamprime) = \muLZ(\gamma_0) \gamprime$, and we have $\gamprime =
\pi(\gamprime)$, so  $d\bigl(\muLZ(\gamma_0 \gamprime) , \pi(H) \bigr) <
C$. But, because $\gamma_0 \gamprime \in \gamma_0\GamPrime$, we have
$\log \|\gamma_0 \gamprime\| \ge d\bigl(\gamma_0, \GamPrime) > R$.
Because $R$ can be made arbitrarily large, while $C$ is fixed, this
contradicts the fact that $\Gamma$ acts properly on $G/H$.

Now consider the general case, where $\GamPrime$ is not assumed to be
contained in~$Z$. For convenience, define $f
\colon \Gamma \to \real^+$ by 
 $f(\gamma) = d\bigl(\muLZ(\gamma) , \pi(H) \bigr)$, and, for $\gamma
\in \GamPrime$, let $\ell(\gamma)$ be the word length of~$\gamma$ with
respect to some (fixed) finite generating set of~$\GamPrime$. To a good
approximation, the argument of the preceding paragraph holds. 
Note that $\ell(\gamma')$ is of order~$f(\gamma_0)$.
 The map $\GamPrime \to \real^+ \colon
 \gamma \mapsto \| \nu (\gamma) \|$
 is bounded above by a polynomial function of~$\ell(\gamma)$ (because
$\nu(\Gamma) \subset N$ consists of unipotent matrices), so 
 $$ \| \nu(\gamprime) \|
 = \bigl( \ell(\gamprime) \bigr)^{O(1)}
 = \bigl( f(\gamma_0) \bigr)^{O(1)} .$$
 Therefore, from Lemma~\ref{linear-change}, we have
 \begin{eqnarray*}
 d \Bigl( \muLZ\bigl(\gamma_0 \pi(\gamprime) \bigr),
\muLZ(\gamma_0 \gamprime) \Bigr)
 &=& O(1) +  \log \|\pi(\gamprime)^{-1} \gamprime \|
 + \log  \| (\gamprime)^{-1} \pi(\gamprime) \| \\
 &=& O(1) + \log \|\nu(\gamprime)\| + \log \|\nu(\gamprime)^{-1}\| \\
 &=& O(1) + O \bigl( \log^+ f(\gamma_0) \bigr) ,
 \end{eqnarray*}
 so
 \begin{eqnarray*}
f(\gamma_0 \gamprime)
 &\le& f\bigl(\gamma_0 \pi(\gamprime) \bigr) + d \Bigl( \muLZ\bigl(\gamma_0
\pi(\gamprime) \bigr), \muLZ(\gamma_0 \gamprime) \Bigr) \\
 &<& C +  \Bigl( O(1) + O\bigl( \log^{+} f(\gamma_0) \bigr) \Bigr)
 < \frac{f(\gamma_0)}{2} ,
 \end{eqnarray*}
  if $f(\gamma_0) > C_1$, where $C_1$ is an appropriate constant that
is independent of $R$, $\gamma_0$, and~$\gamprime$.

 This is the start of an inductive procedure: given $\gamma_0$ with
$d(\gamma_0, \GamPrime) > R$, construct a sequence
$\gamma_0,\gamma_1,\ldots,\gamma_M$ in $\gamma_0 \GamPrime$, such
that $f(\gamma_{n+1}) < f(\gamma_n)/2$, for each~$n$. Terminate
the sequence when $f(\gamma_M) \le C_1$, which, obviously, must
happen after only finitely many steps. However, because $\log
\|\gamma_M\| \ge d(\gamma_M , \GamPrime) = d(\gamma_0,\GamPrime) > R$ is
arbitrarily large, this contradicts the fact that $\Gamma$ acts properly
on $G/H$.
 \end{proof}

\begin{rem} \label{Benoist-nilp}
 We now describe how to eliminate the simplifying assumptions made at the
start of the proof of Theorem~\ref{Benoist-NinH}. We discuss only one
assumption at a time; the general case is handled by employing a combination
of the arguments below.
 \begin{enumerate}

\item \emph{Suppose that $\Gamma$ is not abelian.}
 Because $\Gamma$ is nilpotent, its Zariski closure (if connected) is a
direct product $T \times U \times E$, where $T$~is a split torus, $U$~is
unipotent, and $E$~is a compact torus. Thus, the beginning of the proof
remains valid without essential change, up to (but not including) the
definition of~$\GamPrime$.

Choose a co-compact, discrete, free abelian subgroup~$\bar\Gamma$ of
$A/(H \cap A)$ that is contained in $\pi(\Gamma)$, choose $\gamma_1,
\gamma_2, \ldots, \gamma_n \in \Gamma$, such that 
$\pi(\gamma_1),\pi(\gamma_2),\ldots,\pi(\gamma_n)$ is a basis
of~$\bar\Gamma$, and let $\GamPrime$ be the subgroup
generated by $\gamma_1, \gamma_2, \ldots, \gamma_n$. (We remark that
$\GamPrime$ may not be closed, because $[\GamPrime,\GamPrime]$ need not
be closed.)

If the closure of~$\GamPrime$ is not co-compact in~$\Gamma$, then
essentially no changes are needed in the proof. Thus, let us suppose that
that the closure of~$\GamPrime$ is co-compact in~$\Gamma$. Then there is
no harm in assuming that $\GamPrime$ is dense in~$\Gamma$, so
$\GamPrime$ is not abelian. 

If $k$ is nonarchimedean, then the closure of every finitely generated
unipotent subgroup is compact; thus, by replacing $\Gamma$ with its
projection into $X \cap A$, we may assume that $\Gamma$ is abelian, so the
original proof is valid.

We may now assume that $k$ is archimedean. Because the closure of
$[\Gamma,\Gamma]$ is a nontrivial unipotent group, we know that it is
noncompact. Thus, there is some $\gamma_0 \in [\Gamma,\Gamma]$ with $\|
\gamma_0 \| > R$. 

Choose $\gamprime \in \GamPrime$, such that $d \bigl( \muLZ(\gamma_0)
\pi(\gamprime), \pi(H) \bigr) < C$.
By choosing $\gamprime$ efficiently, we may assume that
$\ell(\gamprime) = O \bigr( \log( \|\gamma_0\|) \bigr)$. This is the
inductive step in the construction of
$\gamma_0,\gamma_1,\ldots,\gamma_M$. Because $\ell(\gamma_n^{-1}
\gamma_{n+1}) = O \bigr( \log( \|\gamma_n^{-1} \gamma_{n+1}\|) \bigr)$,
we see that 
 $$\| \nu(\gamma_0^{-1} \gamma_M) \| = O( M \log \|
\gamma_0\|) = O \bigl( (\log \|
\gamma_0\|)^2 \bigr) \ll \| \gamma_0\| = \| \nu(\gamma_0) \| .$$
 Therefore $\| \gamma_M \| \ge \|
\nu(\gamma_M) \| \approx \| \nu(\gamma_0) \|$ is large. So $\Gamma$ is not
proper on $G/H$.

\item \emph{Do not assume that $H = (H \cap A) \semiprod (H \cap
N)$.}
 There is no harm in assuming that $H \subset AN$ (see~\ref{HcanbeAN} or
\cite[Lem.~4.2.2]{Benoist}). If $H$ is Zariski closed (and Zariski
connected), then, after replacing $H$ by a conjugate subgroup, we have $H =
(H \cap A) \semiprod (H \cap N)$. Thus, the problem only arises when $H$ is
not Zariski closed. In this case, $H$ must be almost connected (and $k$ must
be archimedean).

The proof remains unchanged up until the choice of~$\gamprime$. Instead of
only choosing an element $\gamprime \in \GamPrime$, we also choose an
element $h' \in H$. Namely, let $h_0 = e$, and choose $\gamprime \in
\GamPrime$ and $h' \in H$, such that
 $$d \bigl( \muLZ(\gamma_0) \pi(\gamprime), \muLZ(h_0) \pi(h') \bigr) < C .$$
 This begins the inductive construction of sequences
$\gamma_0,\gamma_1,\ldots,\gamma_M$ in~$\Gamma$ and $h_0,h_1,\ldots,h_M$
in~$H$, such that $d(\gamma_n, e) > R$ for each~$n$, and 
 $$d \bigl( \muLZ(\gamma_M), \muLZ(h_M) \bigr) < C_2 ,$$
 for an appropriate constant~$C_2$ that is independent of $R$. This
contradicts the fact that $\Gamma$ is proper on $G/H$.

 \item \emph{Suppose that $\Gamma \backslash G/H$ is not compact.}
 The compactness was used only to show that
$A/ \closure{\pi(\Gamma H)}$ is compact, and that if $\GamPrime$ is
co-compact in~$\Gamma$, then $H$ contains~$N$. Thus, it suffices to show,
after replacing $H$ by a conjugate, that \pref{Benoist-nilp-A/cpct}~$A/
\closure{\pi \bigl(\Gamma (H \cap AN) \bigr)}$ is compact, and that
\pref{Benoist-nilp-NinH}~if $\GamPrime$ is co-compact in~$\Gamma$, then
$H$ contains a conjugate of~$N$.

 The finite measure on
$\Gamma \backslash G/H$ pushes to a finite measure on the quotient $AN
\backslash G/H$. (We may assume that $\Rad H \subset AN$, so the quotient
$AN \backslash G / H$ is countably separated.) By considering a generic fiber
of this quotient map, we see that we may assume, after replacing $H$ by a
conjugate subgroup, that $\Gamma \backslash AN / (H \cap AN)$ has finite
volume.

\begin{enumerate}

 \item \label{Benoist-nilp-A/cpct}
 Because $A/\closure{\pi\bigl(\Gamma(H
\cap AN)\bigr)}$ is a quotient of $\Gamma \backslash AN / (H \cap AN)$, it
must have finite volume. Therefore, it is compact.

\item \label{Benoist-nilp-NinH}
 Suppose that $\GamPrime$ is co-compact in~$\Gamma$. Almost every fiber of
the quotient map $\Gamma \backslash AN / (H \cap AN) \to A/\pi\bigl(\Gamma(H
\cap AN)\bigr)$ must have finite measure. Thus, we may assume that $\Gamma
\backslash \Gamma N (H \cap AN)/(H \cap AN)$ has finite measure. Because
$\GamPrime$ is co-compact in~$\Gamma$ (and $\Gamma$ is discrete), this
implies that $N (H \cap AN)/(H \cap AN)$ has finite measure, or,
equivalently, that $N/(H \cap N)$ has finite measure. Therefore, $H \cap N =
N$ \see{BDT-solv}, so $N \subset H$, as desired.

 \end{enumerate}
 \end{enumerate}
 \end{rem}

\begin{prop}[{Benoist, cf.~\cite[Prop.~5.1]{Benoist}}] \label{linear-change}
 There is a constant $C>0$ such that, for all $g,h \in G$, we have
 $$ \max \bigl\{ \| \mu(g)^{-1} \mu(gh) \|, \| \mu(g)^{-1} \mu(hg) \|
\bigr\}
 \le C \max \bigl\{ \| h \|, \| h^{-1} \| \bigr\} .$$
 \end{prop}

\end{appendix}

\end{document}